%% file: log.tex
\documentclass[12pt]{article}

\usepackage{graphicx}
\usepackage{latexsym,amsmath,amsfonts,amscd, amsthm}
\usepackage{color}
\usepackage{bm}
\usepackage{tikz}
\usepackage{multirow}
\usetikzlibrary{arrows,backgrounds,snakes}

\newtheoremstyle{plainNoItalics}{}{}{\normalfont}{}{\bfseries}{.}{ }{}

\theoremstyle{plain}
\newtheorem{thm}{Theorem}[section]

\theoremstyle{plainNoItalics}

\newtheorem{defn}[thm]{Definition}
\newtheorem{rem}[thm]{Remark}
\newtheorem{prop}[thm]{Proposition}

\renewcommand{\theequation}{\thesection.\arabic{equation}}

\newcommand{\f}{\frac}

\newcommand{\beq}{\begin{equation}}
\newcommand{\eeq}{\end{equation}}
\newcommand{\beqa}{\begin{eqnarray}}
\newcommand{\eeqa}{\end{eqnarray}}
\newcommand{\bit}{\begin{itemize}}
\newcommand{\eit}{\end{itemize}}
\newcommand{\bedef}{\begin{defn}}
\newcommand{\edefn}{\end{defn}}
\newcommand{\bpro}{\begin{prop}}
\newcommand{\epro}{\end{prop}}

\begin{document}

\baselineskip=1.8pc

%\vspace*{.10in}

%=============  title  =========================

\input{title}

\newpage

\input{intro}
\input{review}
\input{numerical}
\input{conclusion}

\input{acknowledge}

\input{appendix}

\bibliographystyle{siam}
\bibliography{refer}

\end{document}

%% file: title.tex
\begin{center}
{\bf
Runge-Kutta Central Discontinuous Galerkin BGK Method for the Navier-Stokes Equations
%Viscous Flow Equations
}
\end{center}
\vspace{.2in}
\centerline{ Tan Ren\footnote{School of
Aerospace Engineering, Beijing Institute of Technology, Beijing, 100081. E-mail: rentanx@gmail.com},
Jun Hu\footnote{School of
Aerospace Engineering, Beijing Institute of Technology, Beijing, 100081. E-mail: hujun@bit.edu.cn},
Tao Xiong\footnote{Department of Mathematics, University of Houston,
Houston, 77004. E-mail: txiong@math.uh.edu.},
Jing-Mei
Qiu\footnote{Department of Mathematics, University of Houston,
Houston, 77004. E-mail: jingqiu@math.uh.edu.}
}

\bigskip
\centerline{\bf Abstract}
\bigskip
\noindent
In this paper, we propose a Runge-Kutta (RK) central discontinuous Galerkin (CDG) gas-kinetic BGK method for the Navier-Stokes equations. The proposed method is based on the CDG method defined on two sets of overlapping meshes to avoid discontinuous solutions at cell interfaces, 
as well as the gas-kinetic BGK model to evaluate fluxes for both convection and diffusion terms. Redundant representation of the numerical solution in the CDG method offers great convenience in the design of gas-kinetic BGK fluxes. Specifically, the evaluation of fluxes at cell interfaces of 
one set of computational mesh is right inside the cells of the staggered mesh, hence the corresponding particle distribution function for flux evaluation is much simpler than that in existing gas-kinetic BGK methods. As a central scheme, the proposed CDG-BGK has doubled the memory requirement 
as the corresponding DG scheme; on the other hand, {for the convection part,} the CFL time step constraint of the CDG method for numerical stability is 
relatively large compared with that for the DG method. Numerical boundary conditions have to be treated with special care. Numerical examples for 1D and 2D viscous flow simulations are presented to validate the accuracy and robustness of the proposed RK CDG-BGK method.
%Given the continuous flow states on the other inside the meshs, it is very easy to evaluate flux at the mesh interface and the %Gaussian quadrature points.
%The gas distribution function at the mesh interface and the Gaussian quadrature points in one set of mesh are smooth in another set of mesh and can be designed as a simpler form. The flux can be evolved relative easily.
\vfill

\noindent {\bf Keywords:}
Central discontinuous Galerkin method,
gas-kinetic scheme,
BGK model,
Navier-Stokes equations,
Computational fluid dynamics.
\newpage

%% file: intro.tex
\section{Introduction}
\label{sec1}
\setcounter{equation}{0}
\setcounter{figure}{0}
\setcounter{table}{0}

There are two scales in describing compressible flow motions: the kinetic scale via the Boltzmann equation describing the particle distribution function and the hydrodynamic scale via the Euler or Navier-Stokes equations describing macroscopic flow variables such as mass, momentum and energy. In a gas-kinetic representation, all flow variables are moments of particle distribution function; the Euler or Navier-Stokes equations can be derived from taking moments of the Boltzmann equation based on the Chapman-Enskog expansion~\cite{chapman1991mathematical}. In this paper, we are interested in numerically simulating Navier-Stokes equations via the Boltzmann scheme with high order accuracy.

In the past few decades, many computational efforts have been devoted to simulate Euler or Navier-Stokes equations in the field of computational fluid dynamics (CFD). Many of classical numerical methods for Navier-Stokes equations involve solving convection and viscous terms separately on one set of computational grid. For the nonlinear convection term, the design of numerical fluxes at element interfaces is crucial to the success of numerical algorithms. Various approximate Riemann solvers have been proposed to approximate the wave structure of exact Riemann solutions, e.g. Godunov scheme~\cite{godunov1959finite}, the approximate Riemann solvers due to Roe~\cite{roe1981approximate}, Osher~\cite{osher1982upwind},  Harten, Lax and van Leer~\cite{harten1983upstream}, etc. For a summary on this topic, see~\cite{toro2009riemann}. 
For the viscous diffusion term, central difference method is often used. Another approach is to design fluxes at cell interfaces based on integrating particle distribution function in 
the phase space at kinetic scale. The kinetic flux vector splitting method for the Euler equations (KFVS-Euler) based on the collisionless Boltzmann equation is introduced in~\cite{mandal1994kinetic}. When viscous effect is considered, the particle distribution function contains both equilibrium (Maxwellian) and nonequilibrium parts of the gas flow; the inviscid and viscous terms can be treated simultaneously. For example, the kinetic flux vector splitting method for Navier-Stokes equations (KFVS-NS) was developed by introducing the nonequilibrium term in particle distribution obtained by Chapman-Enskog expansion in~\cite{chou1997kinetic}; the gas-kinetic BGK method for Navier-Stokes equations {(BGK-NS)} was introduced in~\cite{xu1994numerical, xu2001gas}.

To improve the performance of numerical schemes, high order schemes are introduced in 80's and underwent great development since then. For example, in the finite volume or finite difference framework, there are the second order MUSCL scheme~\cite{van1979towards}, the essentially non-oscillatory (ENO) scheme~\cite{harten1987uniformly} and the weighted ENO (WENO) scheme~\cite{liu1994weighted, jiang1996efficient}. The discontinuous Galerkin (DG) method, as a class of finite element methods, has been very popular in the CFD community~\cite{cockburn1989tvb2, bassi1997high, cockburn2001runge}. The high-order accuracy of DG is achieved by using high-order polynomial approximations within each element, where more than one degrees of freedom per element are stored and updated. 
The DG method has been well-known for its flexibility, h-p adaptivity, compactness and high parallel efficiency~\cite{cockburn2000development}. 
There have been many work in existing literatures in improving the BGK-NS method to be of high order 
accurate by interpolations or reconstructions such as WENO~\cite{luo2013comparison}, by piecewise parabolic reconstruction of high order BGK fluxes~\cite{li2010high}, and by the DG framework~\cite{xu2004discontinuous, liu2007runge, ni2008dgbgk}. These methods have been  successfully applied in many engineering problems, such as the hypersonic viscous and heat conducting flows~\cite{li2005application, xu2005multidimensional, liao2007gas}, 3D transonic flow ~\cite{may2007improved},  among many others. Comparisons between the schemes with approximate Riemann solvers and the particle distribution functions of the Boltzmann equation are provided in~\cite{luo2013comparison, li2011comparison}.

The central scheme uses staggered meshes to avoid solving Riemann problems at cell interfaces and provides a black box solution to nonlinear hyperbolic conservation laws~\cite{nessyahu1990non, kurganov2000new}.  Liu~\cite{liu2005central} introduced central schemes based on two sets of overlapping meshes. Taking advantages of the redundant representation of the solution on overlapping meshes, approximate Riemann solvers are not needed at cell interfaces, and the high order total variation diminishing (TVD) Runge-Kutta (RK) methods can be directly applied by the method of lines approach. Following similar spirit, central DG (CDG) methods are proposed and developed for hyperbolic equations in~\cite{liu2007central}, and central local DG methods are proposed for diffusion equations in~\cite{liu2011central}.

We propose to couple the CDG framework~\cite{liu2007central} with the BGK-NS method~\cite{xu2001gas} for Navier-Stokes simulations. Compared with the DG BGK methods~\cite{xu2004discontinuous, liu2007runge, ni2008dgbgk},
CDG methods evolve two pieces of approximate solutions defined on two sets of overlapping meshes. Such redundant representation of numerical solution offers great convenience in the design of gas-kinetic BGK fluxes. Specifically, the evaluation of fluxes at cell interfaces of one set of computational mesh is right inside a cell of the staggered mesh (i.e. continuous regions of the solution at the staggered mesh). Hence, the particle distribution function, without involving two different Maxwellian distributions from the left and right states and the corresponding equilibrium state, is much simpler than existing gas-kinetic BGK methods~\cite{xu2001gas, xu2004discontinuous, liu2007runge, ni2008dgbgk}. 
One of the key components that contributes to the success of the gas-kinetic BGK scheme~\cite{xu2001gas} is the {\em exact time evolution of the BGK equation along characteristics}; 
such mechanism, despite its rather complicated formulation, brings the distribution function at cell interfaces to the equilibrium 
state in a very effective way. In the CDG framework, since the distribution function is continuous (at the interior of the other set of solution), such exact evolution is not as crucial. In our scheme, the method of lines approach is adopted; a third-order TVD RK method is used for temporal discretization. As the central scheme, the proposed CDG-BGK has doubled the memory requirement, since two sets of solutions have to be stored and updated simultaneously; on the other hand, {for the convection part,} the CFL time step constraint of the CDG method for numerical stability is relatively large compared with that for the DG method. The numerical boundary conditions have to be treated with special care. For example, a class of DG basis functions that preserves the given boundary condition, in the spirit of~\cite{cockburn2004locally}, 
are proposed for the wall boundary condition.

The paper is organized as follows. In Section 2,  we propose the CDG-BGK method for one and two dimensional problems. The BGK fluxes, as well as numerical boundary conditions are discussed in details. In Section~\ref{secnu}, following the pioneering work of~\cite{xu2001gas, liu2007runge}, extensive numerical results are demonstrated to showcase the effectiveness of the proposed approach. We conclude the paper in Section~\ref{seccon}.

%\bibliographystyle{siam}
%\bibliography{refer}
%
%
%\end{document}

%% file: review.tex
\section{CDG-BGK method for compressible Navier-Stokes equations}
\label{sec2}
\setcounter{equation}{0}
\setcounter{figure}{0}
\setcounter{table}{0}

In this section, we first introduce a 1D BGK model and the corresponding macroscopic conservative Navier-Stokes equations in Section~\ref{1DBGK}. We propose to use the central discontinuous Galerkin (CDG) spatial discretization~\cite{liu2005central} coupled with a third-order total variation diminishing (TVD) Runge-Kutta (RK) temporal discretization~\cite{shu1988efficient} for solving the 1D {Navier-Stokes equations} in Section~\ref{semiCDG}. The BGK type flux for both convection and viscous terms will be described in Section~\ref{BGKfluxes} and extension to {two}-dimensional cases will be presented in Section~\ref{2DCDG}. Finally we discuss the numerical boundary conditions in Section~\ref{BC}.

\subsection{The 1D BGK model and Navier-Stokes equations}
\label{1DBGK}
The integro-differential kinetic Boltzmann equation is commonly used to describe the evolution of the particle distribution function. To avoid the complicated bilinear collisional operator of the Boltzmann equation, a simplified BGK model was proposed by Bhatnagar et al.~\cite{bhatnagar1954model}.
The BGK collisional operator is known to preserve the collisional invariant properties of mass, momentum and energy, as well as the entropy dissipation property.

For a 1D flow, the BGK model can be written as \cite{xu2001gas}
\begin{equation}
f_t+uf_x=\frac{g-f}{\tau},
\label{BGK model}
\end{equation}
where $\tau$ is the particle collision time, $f(x, t, u, \bm{\xi})$ is an unknown function of space variable $x$, time variable $t$, particle velocity $u$ and internal variables $\bm{\xi}$. $\bm{\xi}$ is taking into account to describe the internal motions, such as rotation and vibration \cite{xu1998gas}.
$g(x, t, u, \bm{\xi})$ is the Maxwellian distribution given by
\begin{equation}
g= \rho \left( \frac{\lambda}{\pi} \right)^{\frac{K+1}{2}} e^{-\lambda \left[(u-U)^2+|\bm{\xi}|^2 \right]},
\label{eq_maxw}
\end{equation}
where $\rho$ is the macroscopic density, 
%where \QQ{$R=1/2$ is the nondimensional gas constant,
$U$ is the macroscopic velocity in the $x$ direction, $\lambda$ is related to the gas temperature $T$ by $\lambda=1/T$, $|\bm{\xi}|^2=\xi_1^2+\xi_2^2+\ldots+\xi_K^2$ with $K$ being the total number of degrees of freedom in $\bm{\xi}$. % \QQ{and each of its component is within $(-\infty, \infty)$.} %\QQQ{What is $\xi$?}

The relation between the macroscopic conservative variables and the microscopic distribution function $f$ is
\begin{equation}
\mathbf{W}=(\rho,\rho U, E)^T=\int \bm{\psi} f\, \mathrm{d}\Xi=\int \bm{\psi} g\, \mathrm{d}\Xi,
\label{conserve_variables}
\end{equation}
where $E=\frac{1}{2}\rho U^2+p/(\gamma-1)$ is the total energy, with $p=\f12\rho T=\rho/(2 \lambda)$ to be the pressure, and $\gamma=(K+3)/(K+1)$ is the ratio of specific heats.
\begin{equation}
\bm{\psi}=(\psi_1,\psi_2,\psi_3)^T=\left(1,u,\frac{1}{2}(u^2+|\bm{\xi}|^2)\right)^T,
\end{equation}
and $\mathrm{d}\Xi=\mathrm{d}u\mathrm{d}\bm{\xi}$ is the volume element in the phase space with $\mathrm{d}\bm{\xi}=\mathrm{d}\xi_1\mathrm{d}\xi_2 \ldots \mathrm{d}\xi_K$.
The second equality in equation \eqref{conserve_variables} is due to the fact that the BGK collisional term conserves mass, momentum and energy. In other words, $f$ and $g$ satisfy the conservation constraint
\begin{equation}
\int \bm{\psi} \frac{g-f}{\tau} \mathrm{d}\Xi=0,
\label{compatibility}
\end{equation}
at any point in space and time. By taking moments of $\bm{\psi}$ to the BGK model~(\ref{BGK model}), due to equation~(\ref{compatibility}), we get
\begin{equation}
\int \bm{\psi} f_t  \mathrm{d}\Xi+  \int  u \bm{\psi} f_x \mathrm{d}\Xi=0,
\label{simpns}
\end{equation}
or
\begin{equation}
\mathbf{W}_t+\mathbf{G}_x=0,
\label{1d_NS_simplified}
\end{equation}
where ${\bf W}$ is the vector of the macroscopic conservative variables in equation~(\ref{conserve_variables}). ${\bf G} =\int u \bm{\psi} f \mathrm{d}\Xi$ is the flux function from the kinetic formulation. Specifically,
\begin{equation}
{\bf G} =
\left(\begin{array}{c}
\displaystyle G_{\rho} \\[3mm]
\displaystyle G_m \\[3mm]
\displaystyle G_E
\end{array}\right) = \int u \left(\begin{array}{c} \displaystyle 1 \\[3mm]
\displaystyle  u \\[3mm]
\displaystyle \frac{1}{2} (u^2+|\bm{\xi}|^2) \end{array}\right) f \mathrm{d} \Xi,
\label{numerical flux}
\end{equation}
where $G_\rho$ is the density flux, $G_m$ is the momentum flux, $G_E$ is the energy flux.

%and then
%\begin{equation}
%\int (g_t+u g_x) \bm{\psi} \mathrm{d}\Xi=\tau \int (g_{tt}+2u g_{xt}+u^2 g_{xx}) \bm{\psi} \mathrm{d}\Xi,
%\end{equation}

%The Chapman-Enskog expansion~\cite{chapman1991mathematical} gives
%\QQ{From the BGK model \eqref{BGK model}, we have}
Based on the Chapman-Enskog expansion~\cite{chapman1991mathematical} with
\beq
\label{eq: f_hilbert}
f=g-\tau(g_t+u g_x) + \mathcal{O}(\tau^2),
\eeq
%\QQQ{From eq.~(\ref{simpns}) as well as the form of the Maxwellian distribution~\eqref{eq_maxw}},
from the BGK model~\eqref{BGK model}, the compressible Navier-Stokes equations {on macroscopic variables} can be obtained by omitting $\mathcal{O}(\tau^2)$ terms (for details see~\cite{xu1998gas})
\begin{equation}
\left(\begin{array}{c}
\displaystyle \rho \\[3mm]
\displaystyle \rho U\\[3mm]
\displaystyle E
\end{array}\right)_t + \left(\begin{array}{c} \displaystyle \rho U \\[3mm]
\displaystyle \rho U^2 + p \\[3mm]
\displaystyle U(E+p) \end{array}\right)_x=
\left(\begin{array}{c} \displaystyle 0 \\[3mm]
\displaystyle s_{1x} \\[3mm]
\displaystyle s_{2x} \\[3mm]\end{array}\right)_x,
\label{Navier-Stokes}
\end{equation}
where $s_{1x}=\mu[ \frac{2K}{K+1} U_x ]$, $s_{2x}=\mu [ \frac{K+3}{4} T_x +\frac{2K}{K+1} U U_x ]$
are the viscous terms, $\mu=\tau p$ is the dynamic viscousity coefficient.
%(see the Appendix~\ref{APP1} for detailed description of these physical variables).
For a monatomic gas, $K=2$, $\gamma=5/3$, the above Navier-Stokes equations become,
\begin{equation}
\left(\begin{array}{c}
\displaystyle \rho \\[3mm]
\displaystyle \rho U\\[3mm]
\displaystyle E
\end{array}\right)_t + \left(\begin{array}{c} \displaystyle \rho U \\[3mm]
\displaystyle \rho U^2 + p \\[3mm]
\displaystyle U(E+p) \end{array}\right)_x=
\left(\begin{array}{c} \displaystyle 0 \\[3mm]
\displaystyle \frac{4}{3} \mu U_x \\[3mm]
\displaystyle \frac{5}{4} \mu T_x+\frac{4}{3} \mu U U_x \\[3mm]\end{array}\right)_x.
\label{Navier-Stokes2}
\end{equation}

\begin{rem}
\label{rem_cd}
The kinetic flux function $\mathbf{G}$ in equation~\eqref{numerical flux} encompasses both the convection and diffusion terms in the macroscopic Navier-Stokes system~\eqref{Navier-Stokes} with the approximation to the $f$ function by equation \eqref{eq: f_hilbert}. The convection term in the Navier-Stokes system~\eqref{Navier-Stokes} is due to the contribution from the Maxwellian function $g$ in \eqref{eq: f_hilbert}, while the diffusion term is due to the $\mathcal{O}(\tau)$ term in equation~\eqref{eq: f_hilbert}.
\end{rem}
% in the form of equation~\eqref{eq_maxw}

\subsection{The RK CDG method}
\label{semiCDG}
We propose to use the RK CDG method~\cite{liu2007central} to solve equation~(\ref{1d_NS_simplified}).
The CDG method evolves two sets of approximate solutions defined on overlapping cells.
Compared with the DG method, the CDG method does not need a numerical flux at the cell interface. %old: \QQ{to define the interface values}.
The evaluation of the flux at the interface of one cell is right inside a cell of the other staggered cells. The CDG method uses the flux function of the solution at the staggered cells, which has no ambiguous
values there. This is convenient for us to define the BGK flux in the next subsection. In the following,
we first follow~\cite{liu2007central} to describe the CDG method.

We first consider a 1D domain $[0,L]$ with a partition of $\{x_i\}_{i=1}^{i=N}$. Denote $x_{i+\frac{1}{2}}=\frac{1}{2}(x_i+x_{i+1})$, and
let $I_i=[ x_{i-\frac{1}{2}}, x_{i+\frac{1}{2}} ]$ and $I_{i+\frac{1}{2}}=[ x_i, x_{i+1} ]$ be two sets of overlapping {cells}.
Two discrete spaces associated with the overlapping {cells} $I_i$ and $I_{i+\frac{1}{2}}$ are defined as
\begin{equation}
\begin{split}
& \mathbf{Z}_h=\mathbf{Z}_h^k= \{\mathbf{z}: \mbox{each of its 3 components}\, z|_{I_i} \in P^k(I_i), \forall i \}, \\
& \mathbf{W}_h=\mathbf{W}_h^k= \{\mathbf{z}: \mbox{each of its 3 components}\, z|_{I_{i+\frac{1}{2}}} \in P^k(I_{i+\frac{1}{2}}), \forall i \}, \notag
\end{split}
\end{equation}
where the local space $P^k(I)$ consists of polynomials of degree at most $k$ on $I$.
%\textcolor{red}{the notations are a little bit strange here; as the W is the vector, while $\eta$ ... are %scalars. Same comments for the 2D case.}

The semi-discrete CDG method for solving~\eqref{1d_NS_simplified} is given as follows: find two sets of approximate solutions
$\mathbf{W}_h^Z \in \mathbf{Z}_h$ and
$\mathbf{W}_h^W \in \mathbf{W}_h$, such that for any $\bm{\eta}_h \in \mathbf{Z}_h$, $\bm{\zeta}_h \in \mathbf{W}_h$ and for all $i$,
\begin{equation}
\begin{split}
\frac{\mathrm{d}}{\mathrm{d} t} \int_{I_i} \mathbf{W}_h^Z \bm{\eta}_h\, \mathrm{d}x &= \frac{1}{\Delta \tau^n} \int_{I_i} (\mathbf{W}_h^W-\mathbf{W}_h^Z) \bm{\eta}_h\, \mathrm{d}x
+\int_{I_i} \mathbf{G}(\mathbf{W}_h^W) \frac{d}{dx}\bm{\eta}_h\, \mathrm{d}x\\
& - \mathbf{G}(\mathbf{W}_h^W(x_{i+\frac{1}{2}},t)) \bm{\eta}_h(x_{i+\frac{1}{2}}^{-}) +\mathbf{G}(\mathbf{W}_h^W(x_{i-\frac{1}{2}},t)) \bm{\eta}_h(x_{i-\frac{1}{2}}^{+}),
\end{split}
\label{semi1d_1}
\end{equation}
\begin{equation}
\begin{split}
 \frac{\mathrm{d}}{\mathrm{d} t} \int_{I_{ i+\frac{1}{2}}} \mathbf{W}_h^W \bm{\zeta}_h\, \mathrm{d}x &= \frac{1}{\Delta \tau^n} \int_{I_{ i+\frac{1}{2}}} (\mathbf{W}_h^Z-\mathbf{W}_h^W) \bm{\zeta}_h\, \mathrm{d}x
 +\int_{I_{ i+\frac{1}{2}}} \mathbf{G}(\mathbf{W}_h^Z) \frac{d}{dx}\bm{\zeta}_h\, \mathrm{d}x\\
& - \mathbf{G}(\mathbf{W}_h^Z(x_{i+1},t)) \bm{\zeta}_h(x_{i+1}^{-}) +\mathbf{G}(\mathbf{W}_h^Z(x_i,t)) \bm{\zeta}_h(x_i^{+}),
\end{split}
\label{semi1d_2}
\end{equation}
where $x_i^\pm$ are the right and left limits at the point $x_i$. {Here, the operations for vectors are in the component-wise sense \cite{cockburn1989tvb2}.}
$\Delta \tau^n$ is the maximum time step determined by the CFL condition, the specification of which can be found at the beginning of Section~\ref{secnu}. The first terms on the right side of equations~(\ref{semi1d_1}) and~(\ref{semi1d_2}) are used to remove its $\mathcal{O}(1/\Delta t)$ dependency of numerical dissipation~\cite{liu2007central}.

We focus our discussions on the approximate solution $\mathbf{W}_h^Z$ on {cell} $I_i$. It can be expressed as
 \begin{equation}
\mathbf{W}_h^Z(x,t)=\sum_{l=0}^{k} \mathbf{W}_i^{Z,l}(t)\eta_{i}^l(x), %\quad \mathbf{W}_h^W(x,t)=\sum_{l=0}^{k} \mathbf{W}_i^{W,l}(t),  and $\{\zeta_i^l \}$
\end{equation}
where $\{\eta_i^l\}$ is a basis function of $P^k(I_i)$. For example, the Legendre polynomials are a local orthogonal basis of $P^k(I_i)$. In the 1D case,
$\eta_i^0=1$, $\eta_i^1=( \frac{x-x_i}{\Delta x_i /2} )$, $\eta_i^2=( \frac{x-x_i}{\Delta x_i /2} )^2-\frac{1}{3}$,$\ldots$, where $\Delta x_i=x_{i+1/2}-x_{i-1/2}$.
The approximate solution $\mathbf{W}_h^W$ on {cell} $I_{i+\frac{1}{2}}$ can be defined similarly.
{The fluxes $\mathbf{G}$
at the cell interfaces and the integrals on the right side of equations~(\ref{semi1d_1}) and (\ref{semi1d_2}) are calculated by the gas-kinetic formulation presented in next section.
The integrals on the right side of equation~\eqref{semi1d_1} contain solutions $\mathbf{W}_h^W, \mathbf{W}_h^Z$, which are continuous over two subintervals $[x_{i-\frac12}, x_i]$ and $[x_i, x_{i+\frac12}]$. These integrals are computed by Gaussian quadrature rules on each of the subinterval. Similar comments apply for the integral on the right hand side of equation~\eqref{semi1d_2}.}
If the solutions are discontinuous, the TVB limiter proposed by Cockburn and Shu~\cite{cockburn1989tvb2} will be used to eliminate spurious oscillations and enforce the stability.

%\QQQ{\begin{rem}
%The central DG method evolves two pieces of approximate solutions defined on overlapping cells, and they avoid the use of Riemann solvers which are complicated and costly for
%system of equations~\cite{leveque2002finite}. Unlike the DG method, the central DG method does not need a numerical flux to define the interface values of the solution, so the evaluation of the solution at
%the interface is in the middle of the staggered cell for uniform cells, hence in the continuous region of the solution.
%For the proposed CDG-BGK scheme, it provides great convenience in terms of scheme design and implementation. Please see the next subsection for more discussions.
%\end{rem}
%}

%The semi-discrete equations~(\ref{semi1d_1}) and (\ref{semi1d_2}) %are ordinary differential equations and
%\QQ{can be written in a concise ordinary differential equations (ODE) form $\frac{d}{dt} \mathbf{W}_h=L_h (\mathbf{W}_h,t)$}, $L_h$ is the spatial operator of the right side and $\mathbf{W}_h=(\mathbf{W}^Z_h, \mathbf{W}^W_h)$.
In this paper, we use a method of lines RK method for temporal discretization. Let $\mathbf{W}_h=(\mathbf{W}^Z_h, \mathbf{W}^W_h)$, the third order TVD RK time discretization~\cite{shu1988total} for equations~(\ref{semi1d_1}) and (\ref{semi1d_2}) is the following,
\begin{equation}
\begin{split}
&\mathbf{W}_h^{(1)}=\mathbf{W}_h^n+\Delta t^n L_h(\mathbf{W}_h^n), \\
&\mathbf{W}_h^{(2)}=\frac{3}{4}\mathbf{W}_h^n+\frac{1}{4}\mathbf{W}_h^{(1)}+\frac{1}{4}\Delta t^n L_h(\mathbf{W}_h^{(1)}), \\
&\mathbf{W}_h^{n+1}=\frac{1}{3}\mathbf{W}_h^n+\frac{2}{3}\mathbf{W}_h^{(2)}+\frac{2}{3}\Delta t^n L_h(\mathbf{W}_h^{(2)}),
\end{split}
\label{rk3}
\end{equation}
where $\Delta t^n=\theta \Delta \tau^n$ is the current time step with $\theta \in (0,1]$
and $L_h (\mathbf{W}_h)$ is the spatial operators on the right side of the semi-discrete equations~(\ref{semi1d_1}) and (\ref{semi1d_2}).

\begin{rem}
%For $P^1,P^2,P^3$ solution spaces, CDG method has a larger CFL number than DG method with the same RK time %discretization. For third-order TVD RK method, CFL numbers of CDG method are $0.58, 0.33, 0.22$ for $P^1, %P^2, P^3$ respectively, while CFL numbers of DG method are $0.40, 0.20, 0.13$ for $P^1, P^2, P^3$ %respectively, which can be found in~\cite{liu20082}.
For the third-order TVD RK method~\eqref{rk3}, the CDG method has a little larger CFL number than the DG method for convection part, e.g., the CFL numbers of the CDG method are $0.58, 0.33, 0.22$, while the CFL numbers of the DG method are $0.4, 0.2, 0.13$, for $P^1, P^2, P^3$ respectively~\cite{liu20082}.
\end{rem}

\subsection{The BGK flux}
\label{BGKfluxes}
In this section, we describe how to evaluate the vector $\bf G$ from the BGK model for the CDG method in equations~(\ref{semi1d_1})-(\ref{semi1d_2}). As pointed out in Remark~\ref{rem_cd}, both convection and viscous terms in Navier-Stokes equations~\eqref{Navier-Stokes} comes from the vector ${\bf G}$, the evaluation of which mimics the Chapman-Enskog expansion~\cite{chapman1991mathematical}.
Below, we provide a brief review of the existing algorithm in numerically approximating the BGK flux for the gas-kinetic scheme in the finite volume framework~\cite{xu2001gas}, with some of the ideas originally developed in~\cite{xu1994numerical, xu1998gas}. Further development of such schemes in the DG framework could be found in~\cite{xu2004discontinuous, liu2007runge, ni2008dgbgk}.
We note that it is very difficult to completely describe the BGK gas-kinetic scheme, hence we only outline main steps with intuition below, but refer readers to the original manuscript ~\cite{xu2001gas} for technical details.
%Usually numerical solutions are discontinuous at cell
%interfaces and a numerical flux is needed to define the flux values at the cell interface.

\underline{Gas-kinetic scheme in finite volume framework~\cite{xu2001gas}.} The gas-kinetic BGK scheme updates the macroscopic conservative variables based on integrating equation~\eqref{1d_NS_simplified} from $t^n$ to $t^{n+1}$,
\beq
\overline{\mathbf{W}}_j^{n+1}=\overline{\mathbf{W}}_j^n+\frac{\Delta t^n}{\Delta x_j}(\widehat{\mathbf{G}}_{j-\frac{1}{2}}^n-\widehat{\mathbf{G}}_{j+\frac{1}{2}}^n),
\eeq
{where $\overline{\mathbf{W}}_j^n$ are cell averages of the macroscopic conservative variables, $\Delta t^n=t^{n+1}-t^n$ is current time step, $\Delta x_j=x_{j+\frac{1}{2}}-x_{j-\frac{1}{2}}$}
and %\QQQ{what is $\overline{W}$? explain. the use of $\Delta t^{n+1}$ is strange.}
{\beq
\widehat{\mathbf{G}}_{j+\frac{1}{2}}^n=\frac{1}{\Delta t^n} \int_{0}^{\Delta t^n} \int \bm{\psi} u f(x_{j+\frac{1}{2}},t^n+t,u,\bm{\xi})\, d\Xi dt.
\label{hatG}
\eeq}
{Here $f$ is to be obtained from analytically solving the BGK model~\eqref{BGK model}
%\beq
%f_t + u f_x = \frac{g-f}{\tau},
%\eeq
via characteristically tracing with a source term. Specifically, %the general solution $f$ at a cell interface $x_{j+1/2}$ and time $t^n+t$ is %\QQ{can be obtained from solving the BGK model \eqref{BGK model} as following,} \QQQ{a little repeated?}
\beq
 f(x_{j+1/2},t^n+t,u,\bm{\xi})=\frac{1}{\tau} \int_{0}^{t} g(x^{\prime},t^{\prime},u,\bm{\xi}) e^{-(t^n+t-t^{\prime})/\tau} \,dt^{\prime} + e^{-t/\tau} f_0 (x_{j+1/2}-u t, t^n, u, \bm{\xi}),
 \label{gkf}
\eeq}
where $x^{\prime}=x_{j+1/2}-u(t^n+t-t^{\prime})$ is the particle trajectory, $f_0$ is the distribution function at $t^n$, $g$ is the Maxwellian distribution function.
%over
%$[x_{j-\frac12}, x_{j+\frac12}]\times[t^n, t^{n+1}]$ to be approximated.
Based on the framework outlined above, to update the macroscopic conservative variable $\overline{\mathbf{W}}_j^{n+1}$, we need to specify procedures to approximate $f_0$ and $g$ in equation~\eqref{gkf} and to evaluate the temporal integration in equation~(\ref{hatG}).
%such way of updating $W^{n+1}$ gas-kinetic BGK model.
%\begin{rem}
%work with BGK scheme, the implementation is at the cost of бн
%\end{rem}

%\QQ{In the following, we discuss the \QQ{evaluation of the BGK} flux at \QQ{at a cell interface} $x_{j+1/2}$ \QQ{at some time $t \in [t^n, t^{n+1}]$}}. %For simplicity, we let $x_{j+1/2}=0$.}
{\em{Approximation of $f_0$.}} %~\cite{chou1997kinetic, xu2001gas}
The initial distribution function is approximated based on equation~\eqref{eq: f_hilbert} as well as a first order Taylor expansion in space around $x_{j+\frac12}$. {We use superscripts $l$ and $r$ to indicate the left and right limits of the solutions respectively. In particular,}
\begin{equation}
 f_0(x,t^n,u,\bm{\xi})=\begin{cases}
    g^l \left[1-\tau(a^l u+A^l)+a^l {(x-x_{j+1/2})}\right], &\mbox{if} \quad {x \le x_{j+1/2}},\\
    g^r \left[1-\tau(a^r u+A^r)+a^r {(x-x_{j+1/2})}\right], &\mbox{if} \quad {x > x_{j+1/2}},
   \end{cases}
  \label{gkf0}
\end{equation}
where %$g^l=g^l(x_{j+1/2},t^n,u,\xi)$ and $g^r=g^r(x_{j+1/2},t^n,u,\xi)$,
\beq
\begin{split}
&g^l=g(x^-_{j+1/2},t^n,u,\bm{\xi})=\rho^l \left( \frac{\lambda^l}{\pi} \right)^{\frac{K+1}{2}} e^{-\lambda^l \left[(u-U^l)^2+|\bm{\xi}|^2
\right]},\\
&g^r=g(x^+_{j+1/2},t^n,u,\bm{\xi})=\rho^r \left( \frac{\lambda^r}{\pi} \right)^{\frac{K+1}{2}} e^{-\lambda^r \left[(u-U^r)^2+|\bm{\xi}|^2
\right]},
\end{split}
\label{glgr}
\eeq
are the Maxwellian distribution functions {based on macroscopic variables $(\rho^l, U^l, \lambda^l)$ and $(\rho^r, U^r, \lambda^r)$} at the left limit $x^-_{j+1/2}$ and the right limit $x^+_{j+1/2}$ respectively. $a^l,  A^l$ are related to the spatial and temporal slopes of $g$ at the left limit $x^-_{j+1/2}$, which are simply denoted as $\partial_x g^l$ and $\partial_t g^l$ with $l$ denoting the left limit. Similarly for $\rho^l_x$, $\lambda^l_x$ and $U^l_x$. Specifically, for $a^l$ we have 
\beq
a^l=\frac{\partial_x g^l}{g^l}=\frac{\rho^l_x}{\rho^l}+\frac{K+1}{2 \lambda^l} \lambda^l_x- \lambda^l_x[(u-U^l)^2+|\bm{\xi}|^2]-2 \lambda^l (U^l-u)U^l_x,
\label{al_macro}
\eeq
which is a quadratic function of $u$ and $\bm{\xi}$. In terms of implementation, it was suggested in~\cite{xu2001gas} to express $a^l$
in the following quadratic form of $u$ and $\bm{\xi}$
\beq
%\begin{split}
a^l=a^l_1+a^l_2 u+a^l_3 \frac{1}{2}(u^2+|\bm{\xi}|^2), \quad
%\quad a^r=a^r_1+a^r_2 u+a^r_3 \frac{1}{2}(u^2+|\xi|^2),\\
%A^l=A^l_1+A^l_2 u+A^l_3 \frac{1}{2}(u^2+|\xi|^2),
%\quad A^r=A^r_1+A^r_2 u+A^r_3 \frac{1}{2}(u^2+|\xi|^2),
%\end{split}
\label{alrAlr}
\eeq
with the coefficients $a^l_1$, $a^l_2$, $a^l_3$ to be determined in a similar fashion as equations~(\ref{eq: a})-(\ref{Gamma}) below.
$A^l$ can also be expressed in the form of
\[
A^l=\frac{\partial_t g^l}{g^l}=A^l_1+A^l_2 u+A^l_3 \frac{1}{2}(u^2+|\bm{\xi}|^2),
\]
whose coefficients can be determined by the compatibility condition as in equations~(\ref{eq: compat})-(\ref{solve_A}) below. Similar notations and comments apply to $a^r,  A^r$. More details can be found in~\cite{xu2001gas}.
%And $\int g^l a^l \psi d\Xi=\frac{\overline{\mathbf{W}}_{j+1/2}^n-\overline{\mathbf{W}}_j^n}{\Delta x/2}$, then $a^l_1, a^l_2, a^l_3$ can be obtained. Similarly, $a^r$ can be determined.
%are  which are corresponding to the
%spatial slopes of the conservative variables obtained by MUSCL reconstruction method, see~\cite{chou1997kinetic, xu2001gas} for details.
%\QQ{For this case see~\cite{chou1997kinetic, xu2001gas} for details, e.g. the spatial slopes are obtained by a MUSCL reconstruction. We will specify how to get the spatial and temporal slopes only for one of the following cases.} % at the end of this section. }

{\em{Approximation of the equilibrium state $g$.}}
{In order to get the Maxwellian function $g$ at $(x^\prime, t^\prime)$ in equation~\eqref{gkf}, let $g$ be approximated by the following} $g_f$ function based on a Taylor expansion of the Maxwellian function around $(x_{j+\frac12}, t^n)$ both in space and in time.
\beq
 g_f(x,t,u,\bm{\xi})=\begin{cases}
    g_0 \left[1+\bar{a}^l (x-x_{j+1/2})+\bar{A} (t-t^n) \right], &\mbox{if} \quad {x \le x_{j+1/2}},\\
    g_0 \left[1+\bar{a}^r (x-x_{j+1/2})+\bar{A} (t-t^n) \right], &\mbox{if} \quad {x > x_{j+1/2}},
   \end{cases}
\label{gkg0}
\eeq
where $g_0$ is the Maxwellian distribution function at $(x_{j+\frac12}, t^n)$ with
\beq
g_0=g_0(x_{j+1/2},t^n,u,\bm{\xi})=\rho_0 \left( \frac{\lambda_0}{\pi} \right)^{\frac{K+1}{2}} e^{-\lambda_0 \left[(u-U_0)^2+|\bm{\xi}|^2
\right]}.
\label{maxg0}
\eeq
{Here the macroscopic variables $\left(\rho_0, (\rho U)_0, E_0 \right)$ at the cell interface $x_{j+1/2}$ and time $t^n$ are obtained by
%by taking the limit $t \to 0 $ in equation~\eqref{gkf} and substituting its solution into eq.~\eqref{compatibility},
\beq
\left( \rho_0, (\rho U)_0, E_0 \right)^T=\int_{u \ge 0} \int g^l \bm{\psi} d\Xi+ \int_{u < 0} \int g^r \bm{\psi} d\Xi,
\eeq
where $g^l, g^r$ are specified in equation~\eqref{glgr}.}
It can be observed that $g_f$ has discontinuous spatial slope and continuous temporal slope around $(x_{j+\frac12}, t^n)$:
$\bar{a}^l,\bar{a}^r$ are related to the spatial slopes of $g_0$ from the left and right sides of $x_{j+1/2}$ and $\bar{A}$ is the temporal slope of $g_0$.
%\QQQ{$\bar{a}^l,\bar{a}^r$ can be determined by a similar process as that for $a^l, a^r$.}
They can be determined by a similar process as that for $a^l, a^r$ and $A^l, A^r$.
The only difference is that we use the `after-collision' Maxwellian $g_0$ here, while we use the `before-collision' equilibrium state $g^l$ and $g^r$ previously.
%\QQQ{$\bar{A}$ can be obtained by the integration of equation~\eqref{compatibility} at $x_{j+1/2}$ over the whole time step $\Delta t^n$, where  $f$ is from equation~\eqref{gkf} and $g_f$ is from equation~\eqref{gkg0}. The details can be found in~\cite{xu2001gas}.}

{\em{Analytical evaluation of the integration in equation~\eqref{hatG}.}}
{After $f$ is determined, we can get
\beq
\widehat{\mathbf{G}}_{j+\frac{1}{2}}^n
=\frac{1}{\Delta t^n} \int_{0}^{\Delta t^n} \left[ \langle \bm{\psi} u f(x_{j+\frac{1}{2}},t^n+t,u,\bm{\xi}) \rangle_{u \ge 0}+ \langle \bm{\psi} u f(x_{j+\frac{1}{2}},t^n+t,u,\bm{\xi}) \rangle_{u<0} \right]\, dt. \notag
\eeq
where %he integration is separated by $u\ge0$ and $u<0$. Specifically,
%\beq
%\langle \cdot \rangle_{u \ge 0}=\int_{u \ge 0} \cdot du d\bm{\xi},
%\eeq
\beq
\langle \bm{\psi} u f(x_{j+\frac{1}{2}},t^n+t,u,\bm{\xi}) \rangle_{u \ge 0}=\int_{u \ge 0} \int \bm{\psi} u f(x_{j+\frac{1}{2}},t^n+t,u,\bm{\xi})\, d\Xi. \notag
\eeq
Similarly for the $u<0$ integral. The details of the moments evaluation of such integration can be found in Appendix A of~\cite{xu2001gas}. The temporal integration in equation~\eqref{hatG} can also be performed analytically.
}
%Due to different form of the $f_0$ (eq. \eqref{gkf0}) and $g_f$ function (eq. \eqref{gkg0}) to the left and right of the cell boundary, the integration in the velocity space is done separately for $u<0$ and $u\ge0$.
%Analytical form of the solutions using the \QQ{exponential} function and complementary error functions \QQ{$erfc$} are available in the literature.
%For example, the details of such integration can be found in \QQQ{Appendix A of}~\cite{xu2001gas}.

%Finally, we would like to remark that the gas-kinetic scheme~\cite{xu2001gas} uses the kinetic BGK model in deriving and evaluating the flux functions for both convection and viscous terms in Navier-Stokes equations. However, in the implementation, after analytical evaluation of the integration in the phase space, only macroscopic quantities such as density, momentum and energy are involved. Hence, the computational cost is at the level of macroscopic models, which is less computationally expensive than the kinetic simulations due to additional phase space dimensions.

{\underline{Proposed CDG-BGK method.}} We propose to use the CDG method for spatial discretization with the BGK flux for convection and diffusion terms in the Navier-Stokes system. In the CDG framework, two pieces of approximate solutions defined on overlapping cells are evolved.
Since the cell interface on one set of cell is right inside a cell of the staggered cells, the solution at the staggered cell is continuous for flux evaluation without ambiguity, and we have
\begin{equation}
g^l = g^r = g_0, \quad a^l=a^r=a, \quad A^l=A^r=A,
\label{smodis}
\end{equation}
in equation~(\ref{gkf0}).
In other words, the before-collision and after-collision equilibrium states are the same. Due to this fact, the use of exact time evolution formula~\eqref{gkf} is not as crucial. Therefore, we propose to use the third-order explicit TVD RK method for time evolution by the method of lines approach. Below we describe in details several numerical approximations needed for the proposed CDG-BGK method. They are, in some sense, special cases for the BGK scheme reviewed above.

{\em{Approximation of the distribution function $f$}.}
The particle distribution function $f$ at any spatial location {$x$} and time $t$ is approximated by
\begin{equation}
 f({x},t,u,\bm{\xi})=g \left[ 1-\tau (au+A) \right].
 \label{ourdis}
\end{equation}
where $g=g({x},t,u,\bm{\xi})$ is the Maxwellian distribution~\eqref{eq_maxw} at time $t$. Similar to equations \eqref{al_macro} and \eqref{alrAlr}, $a$ is related to the spatial slope of $g$ in the form of
\begin{equation}
\label{eq: a}
\frac{\partial_x g}{g}\doteq a=a_1+a_2u+a_3 \frac12(u^2+|\bm{\xi}|^2).
\end{equation}
The coefficients $a_1$, $a_2$ and $a_3$ are determined by taking the spatial derivatives on the components of $\mathbf{W}$ in equation~\eqref{conserve_variables},
%It has a unique correspondence with the slopes of the macroscopic conservative variables~{~\cite{xu2001gas}},
\begin{equation}
\int a g\, \mathrm{d}\Xi= \frac{\partial \rho}{\partial x}, \quad
\int aug\,  \mathrm{d}\Xi=\frac{\partial (\rho U)}{\partial x}, \quad
\int a \frac{1}{2} (u^2+|\bm{\xi}|^2)g\, \mathrm{d}\Xi=\frac{\partial E}{\partial x},
\label{slopecon}
\end{equation}
where the slopes of the macroscopic conservative variables can be obtained by directly taking derivatives of the CDG polynomials.
Equation~(\ref{slopecon}) can be rewritten in a matrix-vector form as
\begin{equation}
\Gamma (a_1, a_2, a_3)^T =\frac{1}{\rho} \left( \frac{\partial \rho}{\partial x}, \frac{\partial {(\rho U)}}{\partial x}, \frac{\partial E}{\partial x} \right)^T,
\label{solve_a}
\end{equation}
where
\begin{equation}
\left(\Gamma_{\alpha \beta}\right)=\left(\int g \psi_\alpha \psi_\beta\, \mathrm{d}\Xi /\rho \right)=
 \left(\begin{array}{ccc}
\displaystyle 1 & U &\Phi_1\\[3mm]
\displaystyle U & U^2+\frac{1}{2 \lambda} & \Phi_2 \\[3mm]
\displaystyle \Phi_1 & \Phi_2 & \Phi_3
\end{array}\right),\quad \alpha, \beta =1,2,3,
\label{Gamma}
\end{equation}
with
\begin{equation}
\begin{split}
 & \Phi_1=\frac{1}{2} \left(U^2+\frac{K+1}{2 \lambda}\right), \quad
 \Phi_2=\frac{1}{2} \left(U^3+\frac{(K+3)U}{2 \lambda}\right), \quad \\
 & \Phi_3=\frac{1}{4} \left(U^4+\frac{K^2+4 K+3}{4 \lambda^2}+\frac{(2 K+6) U^2}{2 \lambda} \right). \notag
\end{split}
\label{matrixphi}
\end{equation}
$A$ is the temporal slope of $g$ with the following form,
\begin{equation}
\frac{\partial_t g}{g} \doteq A=A_1+A_2u+A_3 \frac12(u^2+|\bm{\xi}|^2),
\label{eq: A}
\end{equation}
where $A_1$, $A_2$, $A_3$ are uniquely determined by the compatibility condition
\beq
\label{eq: compat}
\int (au+A)\bm{\psi} g \, \mathrm{d}\Xi=0,
\eeq
that is
\begin{equation}
\Gamma (A_1, A_2, A_3)^T=-\frac{1}{\rho} \int au \bm{\psi} g \, \mathrm{d}\Xi,
\label{solve_A}
\end{equation}
where $\Gamma$ is the same as equation~\eqref{Gamma}.
The matrix $\Gamma$ is symmetric and can be efficiently inverted to determine the components of $a$ in equation~\eqref{eq: a} and $A$ in equation~\eqref{eq: A}.

{\em{Integration in the phase space to obtain $\mathbf{G}$.}}
After $a, A$ are determined, we can get the fluxes in equations~(\ref{semi1d_1}) and (\ref{semi1d_2}) by taking the moments of $u\bm{\psi}$ to the distribution function $f$ given by equation~\eqref{ourdis},
\beq
\mathbf{G}=\int u \bm{\psi} f d\Xi=\rho \left[ \langle u\bm{\psi} \rangle-\tau \langle au^2\bm{\psi} \rangle- \tau \langle Au \bm{\psi} \rangle \right].
\label{Gdeterm}
\eeq
The moments evaluation of the Maxwellian distribution function are provided in Appendix~\ref{APP2}.
\begin{rem}
Although the derivation of the numerical flux functions $\mathbf{G}$ comes from the kinetic BGK formulation, the integration in the phase space is done analytically as in Appendix~\ref{APP2}.
The actual implementation is at the level of macroscopic variables. Hence, the computational cost is on the same scale of other existing Navier-Stokes solvers. The same comments apply to
the gas-kinetic BGK scheme in \cite{xu2001gas}.
\end{rem}

\begin{rem}
The BGK model corresponds to {a} unit Prandtl number $Pr$. For a variable Prandtl number, we modify the energy flux by~\cite{xu2001gas} % subtracting the heat flux and adding another one with a variable Prandtl number~\cite{xu2001gas}.
\begin{equation}
 G_E^{new}=G_E+(\frac{1}{Pr}-1)q.
 \label{pr_modify}
\end{equation}
The heat flux $q$ can be evaluated precisely,
\begin{equation}
 q=\frac{1}{2} \int \left(u-U \right) \left( \left(u-U \right)^2+\bm{\xi}^2 \right) f \mathrm{d}\Xi
 =-\tau \int g \left(u-U \right) \left(\psi_3-\psi_2 U+\frac{1}{2} U^2 \right) \left(au+A \right) \mathrm{d}\Xi.
 \end{equation}
%This formula is used in all numerical examples except laminar boundary layer case in Section~\ref{secnu}.
\end{rem}

\begin{rem}
In equation~(\ref{slopecon}), the spatial derivative is directly taken on the CDG polynomials.
It would lead to a $k^{th}$ order scheme with $P^k$ polynomial space for viscous flow simulations {(not the optimal $(k+1)^{th}$ order)}.
For convection dominated problems, the accuracy would still be $(k+1)^{th}$ order.
Such fact {is verified} in our numerical results in Section~\ref{secnu}.
%is being verified
\end{rem}

\subsection{Extension to {two-}dimensional cases}
\label{2DCDG}
In this subsection, we extend the proposed CDG-BGK method in previous subsections to {two-}dimensional cases. The 2D BGK model can be written as
\begin{equation}
f_t+uf_x+vf_y=\frac{g-f}{\tau},
\label{2d_bgk}
\end{equation}
{where $f(x,y,t,u,v,\bm{\xi})$ is the particle distribution function of space variables $x, y$, time variable $t$, particle velocities $u, v$ and internal variables $\bm{\xi}$. While $g(x,y,t,u,v,\bm{\xi})$ is the Maxwellian distribution given by}
\begin{equation}
g= \rho \left( \frac{\lambda}{\pi} \right)^{\frac{K+2}{2}} e^{-\lambda[(u-U)^2+(v-V)^2+|\bm{\xi}|^2]},
\label{2dmaxw}
\end{equation}
{where $\rho$ is the macroscopic density, $\lambda=1/(2RT)$, $U,V$ are the macroscopic velocities in $x, y$ directions,
$|\bm{\xi}|^2=\xi_1^2+\xi_2^2+\ldots+\xi_K^2$ with $K$ being the total number of degrees of freedom in $\bm{\xi}$.} %$\bm{\xi}$ and $K$ have the same meanings as the 1D case.
Based on the conservative constraint condition~\eqref{compatibility}, taking the moments of $\bm{\psi}=(\psi_1, \psi_2, \psi_3, \psi_4)^T=(1,u,v,\frac{1}{2}(u^2+v^2+|\bm{\xi}|^2))^T$ to equation~\eqref{2d_bgk},
we can get the following system of macroscopic conservative equation
\begin{equation}
\int \bm{\psi} f_t\, \mathrm{d}\Xi+ \int u \bm{\psi} f_x \, \mathrm{d}\Xi+ \int v \bm{\psi} f_y\, \mathrm{d}\Xi=0,
\label{2dns1}
\end{equation}
or
\begin{equation}
\mathbf{W_t}+\mathbf{G_x}+\mathbf{H_y}=0.
\label{2dns}
\end{equation}
%corresponding to the 2D Navier-Stokes equations when $f$ is approximate to the first order of $\tau$.
Here ${\bf W}=(\rho, \rho U, \rho V, E)^T$ is the vector of macroscopic conservative variables,
${\bf G} =\int u \bm{\psi} f \mathrm{d}\Xi$ and ${\bf H} =\int v \bm{\psi} f \mathrm{d}\Xi$ are {the flux functions in $x, y$ directions respectively.}
%consists of both convection and viscous terms in $x,y$ directions respectively.
The Chapman-Enskog expansion with $f=g-\tau (g_t+ug_x+vg_y)+ \mathcal{O}(\tau^2)$, and from equation~\eqref{2dns1} gives
a 2D compressible Navier-Stokes equations \cite{xu1998gas},
\begin{equation}
\left(\begin{array}{c}
\displaystyle \rho \\[3mm]
\displaystyle \rho U\\[3mm]
\displaystyle \rho V\\[3mm]
\displaystyle E
\end{array}\right)_t +
\left(\begin{array}{c}
\displaystyle \rho U \\[3mm]
\displaystyle \rho U^2 + p \\[3mm]
\displaystyle \rho UV \\[3mm]
\displaystyle U(E+p) \end{array}\right)_x+
\left(\begin{array}{c}
\displaystyle \rho V \\[3mm]
\displaystyle \rho UV \\[3mm]
\displaystyle \rho V^2+p \\[3mm]
\displaystyle V(E+p) \end{array}\right)_y=
\left(\begin{array}{c} \displaystyle 0 \\[3mm]
\displaystyle s_{1x} \\[3mm]
\displaystyle s_{2x} \\[3mm]
\displaystyle s_{3x} \end{array}\right)_x+
\left(\begin{array}{c} \displaystyle 0 \\[3mm]
\displaystyle s_{1y} \\[3mm]
\displaystyle s_{2y} \\[3mm]
\displaystyle s_{3y} \end{array}\right)_y,
\label{2dNavier-Stokes}
\end{equation}
where
\begin{equation}
\begin{split}
 & S_{1x}=\mu \left[2U_x-\frac{2}{K+2} (U_x+V_y) \right],\quad S_{1y}=\mu (U_y+V_x),\\
 & S_{2x}=\mu(V_x+U_y), \quad S_{2y}=\mu \left[ 2V_y-\frac{2}{K+2}(U_x+V_y) \right],\\
 & S_{3x}=\mu \left[ 2UU_x+V(V_x+U_y)-\frac{2}{K+2} U(U_x+V_y) +\frac{K+4}{4} T_x \right],\\
 & S_{3y}=\mu \left[ U(U_y+V_x)+2VV_y-\frac{2}{K+2} V(U_x+V_y) +\frac{K+4}{4} T_y \right]. \notag
 \end{split}
\end{equation}
The total energy $E=\frac{1}{2} \rho (U^2+V^2)+p/(\gamma-1)$ with the pressure $p=\rho/(2 \lambda)$ and $\gamma=(K+4)/(K+2)$. For monatomic gas $K=1$ and $\gamma=5/3$, while for diatomic gas $K=3$ and $\gamma=7/5$. Similarly, the kinetic flux functions $\mathbf{G, H}$ in equation~\eqref{2dns} are represented by both convection and diffusion terms
 in the 2D macroscopic Navier-Stokes equations~\eqref{2dNavier-Stokes} in $x, y$ directions respectively.

\begin{figure}[ht]
\begin{center}
\includegraphics[width=3.in]{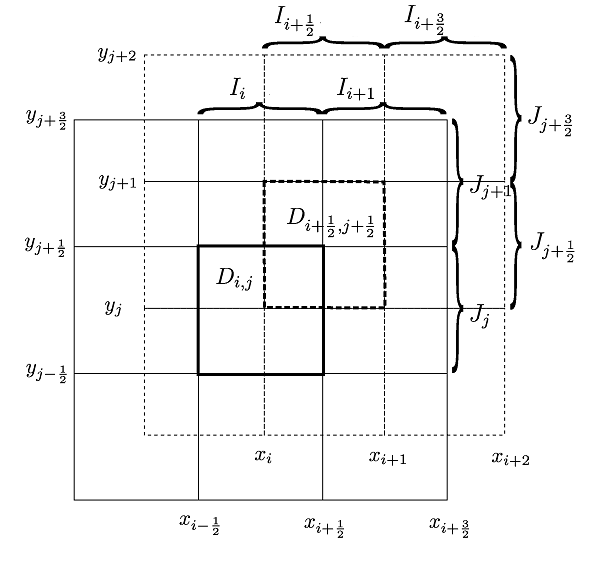}
\end{center} \caption{2D overlapping meshes.}
\label{2doverlapmesh}
\end{figure}

We consider the following numerical discretization of a 2D rectangular domain $\Omega=[0,L_x] \times [0,L_y]$. Let $\{x_i\}_{i=1}^{i=N_x}$ and $\{y_j\}_{j=1}^{j=N_y}$ be partitions of $[0,L_x]$ and $[0,L_y]$ respectively, with
$x_{i+\frac{1}{2}}=\frac{1}{2}(x_i+x_{i+1})$, $y_{j+\frac{1}{2}}=\frac{1}{2}(y_j+y_{j+1})$. Let $I_i=[x_{i-\frac{1}{2}}, x_{i+\frac{1}{2}}]$, $J_j=[y_{j-\frac{1}{2}}, y_{j+\frac{1}{2}}]$,
$I_{i+\frac{1}{2}}=[x_{i}, x_{i+1}]$, $J_{j+\frac{1}{2}}=[y_{j}, y_{j+1}]$.
Denote $\{D_{i,j}\}_{i,j}$ and $\{D_{ i+\frac{1}{2}, j+\frac{1}{2} }\}_{i,j}$ be two sets of overlapping meshes for $\Omega$, with $D_{i,j}=I_i \times J_j$
and $D_{ i+\frac{1}{2}, j+\frac{1}{2} }=I_{i+\frac{1}{2}} \times J_{j+\frac{1}{2}}$, see Fig.~\ref{2doverlapmesh}.
Two discrete spaces associated with the overlapping {cells} $\{D_{i,j}\}_{i,j}$ and $\{D_{ i+\frac{1}{2}, j+\frac{1}{2} }\}_{i,j}$ are defined as
\begin{equation}
\begin{split}
& \mathbf{Z}_h=\mathbf{Z}_h^k= \{\mathbf{z}:\mbox{each of its {4} components}\, z|_{D_{i,j}} \in P^k(D_{i,j}), \forall i,j \}, \\
& \mathbf{W}_h=\mathbf{W}_h^k= \{\mathbf{z}:\mbox{each of its {4} components}\, z|_{D_{i+\frac{1}{2},j+\frac{1}{2}}} \in P^k(D_{i+\frac{1}{2},j+\frac{1}{2}}), \forall i,j \}, \notag
\end{split}
\end{equation}
where the local space $P^k(D)$ consists of polynomials of degree at most $k$ on $D$.
Similar to the 1D case, the semi-discrete CDG scheme for solving equation~\eqref{2dns} is given as follows: find two sets of approximate solutions $\mathbf{W}_h^Z \in \mathbf{Z}_h$ and $\mathbf{W}_h^W \in \mathbf{W}_h$, such that for any $\bm{\eta}_h \in \mathbf{Z}_h$, $\bm{\zeta}_h \in \mathbf{W}_h$ and for all $i$ and $j$,
\begin{equation}
\begin{split}
\frac{\mathrm{d}}{\mathrm{d} t} \int_{D_{i,j}} \mathbf{W}_h^Z \bm{\eta}_h\, \mathrm{d}x \mathrm{d}y&= \frac{1}{\Delta \tau^n} \int_{D_{i,j}} (\mathbf{W}_h^W-\mathbf{W}_h^Z) \bm{\eta}_h\, \mathrm{d}x \mathrm{d}y
+\int_{D_{i,j}} (\mathbf{G}(\mathbf{W}_h^W) \frac{d}{dx}\bm{\eta}_h+\mathbf{H}(\mathbf{W}_h^W) \frac{d}{dy}\bm{\eta}_h)\, \mathrm{d}x  \mathrm{d}y\\
& - \int_{y_{j-1/2}}^{y_{j+1/2}} (\mathbf{G}(\mathbf{W}_h^W(x_{i+\frac{1}{2}},y,t)) \bm{\eta}_h(x_{i+\frac{1}{2}}^{-},y) -\mathbf{G}(\mathbf{W}_h^W(x_{i-\frac{1}{2}},y,t)) \bm{\eta}_h(x_{i-\frac{1}{2}}^{+},y)) \mathrm{d}y \\
& - \int_{x_{i-1/2}}^{x_{i+1/2}} (\mathbf{H}(\mathbf{W}_h^W(x,y_{j+\frac{1}{2}},t)) \bm{\eta}_h(x,y_{j+\frac{1}{2}}^{-}) -\mathbf{H}(\mathbf{W}_h^W(x,y_{j-\frac{1}{2}},t)) \bm{\eta}_h(x,y_{j-\frac{1}{2}}^{+})) \mathrm{d}x, \\
\end{split}
\label{semi2d_1}
\end{equation}
\begin{equation}
\begin{split}
 \frac{\mathrm{d}}{\mathrm{d} t} \int_{D_{ i+\frac{1}{2}, j+\frac{1}{2} }} \mathbf{W}_h^W \bm{\zeta}_h\, \mathrm{d}x \mathrm{d}y&= \frac{1}{\Delta \tau^n} \int_{D_{ i+\frac{1}{2}, j+\frac{1}{2} }} (\mathbf{W}_h^Z-\mathbf{W}_h^W) \bm{\zeta}_h\, \mathrm{d}x \mathrm{d}y\\
&
 +\int_{D_{ i+\frac{1}{2}, j+\frac{1}{2} }} (\mathbf{G}(\mathbf{W}_h^Z) \frac{d}{dx}\bm{\zeta}_h+\mathbf{H}(\mathbf{W}_h^Z) \frac{d}{dy}\bm{\zeta}_h)\, \mathrm{d}x  \mathrm{d}y \\
& - \int_{y_j}^{y_{j+1}} (\mathbf{G}(\mathbf{W}_h^Z(x_{i+1},y,t)) \bm{\zeta}_h(x_{i+1}^{-},y) -\mathbf{G}(\mathbf{W}_h^Z(x_i,y,t)) \bm{\zeta}_h(x_i^{+},y)) \mathrm{d}y \\
& - \int_{x_i}^{x_{i+1}} (\mathbf{H}(\mathbf{W}_h^Z(x,y_{j+1},t)) \bm{\zeta}_h(x,y_{j+1}^{-}) -\mathbf{H}(\mathbf{W}_h^Z(x,y_j,t)) \bm{\zeta}_h(x,y_j^{+})) \mathrm{d}x. \\
\end{split}
\label{semi2d_2}
\end{equation}
Here vectors operations are component-wise operations.
The semi-discrete scheme of equations~(\ref{semi2d_1})-(\ref{semi2d_2}) will be evolved in time by the
third-order TVD RK time method~\eqref{rk3}.
The approximate solution $\mathbf{W}_h^Z$ on the element
$D_{i,j}$ can be expressed as
\begin{equation}
\mathbf{W}_h^Z(x,y,t)=\sum_{l=0}^{2k+1} \mathbf{W}_{i,j}^{Z,l}(t)\eta_{i,j}^l(x,y),  \quad \mbox{for} \quad x,y \in D_{i,j}.
\end{equation}
The 2D Legendre polynomials $\eta_{i,j}^l$ are taken as a local orthogonal basis on $D_{i,j}$,
\begin{equation}
\begin{split}
&\eta_{i,j}^0(x,y)=1,\quad \eta_{i,j}^1(x,y)=\frac{(x-x_i)}{\Delta x_i /2},\quad \eta_{i,j}^2(x,y)=\frac{(y-y_j)}{\Delta y_j /2},
\quad \eta_{i,j}^3(x,y)=\eta_{i,j}^1(x,y) \eta_{i,j}^2(x,y),\\ &\eta_{i,j}^4(x,y)=(\eta_{i,j}^1(x,y)) ^2-\frac{1}{3},\quad \eta_{i,j}^5(x,y)=(\eta_{i,j}^2(x,y)) ^2-\frac{1}{3}, \ldots, \notag
\end{split}
\end{equation}
where $\Delta x_i=x_{i+\frac12}-x_{i-\frac12}$, $\Delta y_j=y_{j+\frac12}-y_{j-\frac12}$. The approximate solution $\mathbf{W}_h^W$ can be obtained similarly.
The fluxes $\mathbf{G}(x,y,t)$ and
$\mathbf{H}(x,y,t)$ on the right side of equations~(\ref{semi2d_1}) and (\ref{semi2d_2}) are calculated by the gas-kinetic formulation presented below.
Note that in 2D case, the first two integrals on the right side of equation~\eqref{semi2d_1} contain four pieces of solutions on four subcells, each of which
 %corresponds to four parts of the solution the flux functions $\mathbf{G}(x,y,t)$ and $\mathbf{H}(x,y,t)$ in each of the volume integral terms corresponds to four parts of the solutions $\mathbf{W}_h^W$ and $\mathbf{W}_h^Z$ respectively.
%and each part
is calculated by a 2D Gaussian quadrature rule. The rest two integrals on the right side of equation~\eqref{semi2d_1} are associated with two pieces of solutions, each of which is determined by a 1D Gaussian quadrature rule. Similar comments apply for the integrals on the right side of equation~\eqref{semi2d_2}.

In the following, we follow the spirit of 1D BGK scheme to propose the 2D strategy for computing the fluxes $\mathbf{G}$ and $\mathbf{H}$.
%\QQ{\underline{Proposed 2D CDG-BGK method.}}
%\QQ{\em{Approximation of the distribution function $f$}.}
The distribution function $f$ in equation~\eqref{2d_bgk} can be expressed as
\begin{equation}
f(x,y,t^n,u,v,\bm{\xi})=g \left[1-\tau (au+bv+A) \right],
\label{gasdis2d}
\end{equation}
up to the first order of $\tau$,
where $g=g(x,y,t^n,u,v,\bm{\xi})$ is the Maxwellian distribution function~\eqref{2dmaxw} associated with the macroscopic variables $(\rho, U, V, \lambda)$ at $(x, y)$ on the dual mesh.
Similar to the 1D case, $a,b$ are related to the spatial slopes of $g$ in $x, y$ directions respectively, and are
%\QQ{Specifically,
%\beq
%a=\partial_x g/g=\partial_x(ln %g)=\frac{\rho_x}{\rho}+\frac{K+2}{2\lambda}\lambda_x-\lambda_x[(u-U)^2+(v-V)^2+|\xi|^2]-2\lambda[(U-u)U_x+(V- %v)V_x], \notag
%\eeq
%which is in a quadratic function of $u, v, \xi$. In terms of implementation,}
%they were suggested in the form of \QQ{~\cite{xu2001gas}}
%\RR{Similar to the 1D case, for implementation, $a$ and $b$ are
taken to be in the form of
\begin{equation}
\frac{\partial_x g}{g} \doteq a=a_1+a_2 u+a_3 v+a_4 \frac12 (u^2+v^2+|\bm{\xi}|^2), \quad
\frac{\partial_y g}{g} \doteq b=b_1+b_2 u+b_3 v+b_4 \frac12 (u^2+v^2+|\bm{\xi}|^2). \notag
\end{equation}
The components of $a$ and $b$ can be uniquely determined from the partial derivatives of the macroscopic conservative variables with respect to $x,y$ %\QQ{~\cite{xu2001gas}},
\begin{equation}
\begin{split}
\int a g\, \mathrm{d}\Xi=\frac{\partial \rho}{\partial x} , \quad &\int b g\, \mathrm{d}\Xi=\frac{\partial \rho}{\partial y} ,\\
\int aug\,  \mathrm{d}\Xi=\frac{\partial (\rho U)}{\partial x}, \quad &\int bug\,  \mathrm{d}\Xi=\frac{\partial (\rho U)}{\partial y},\\
\int avg\,  \mathrm{d}\Xi=\frac{\partial (\rho V)}{\partial x}, \quad &\int bvg\,  \mathrm{d}\Xi=\frac{\partial (\rho V)}{\partial y},\\
\int a \frac{1}{2} (u^2+v^2+|\bm{\xi}|^2)g\,  \mathrm{d}\Xi=\frac{\partial E}{\partial x}, \quad &\int b \frac{1}{2} (u^2+v^2+|\bm{\xi}|^2)g\,  \mathrm{d}\Xi=\frac{\partial E}{\partial y}.
\end{split}
\label{2dslopecon}
\end{equation}
The above equations can be written in a matrix-vector form as
%We focus on the determination process of $a$, $b$ can be obtained similarly. By integrating the left side of equation~(\ref{2dslopecon}), it can be written as
\begin{equation}
\Gamma (a_1, a_2, a_3, a_4)^T=\frac{1}{\rho} \left(\frac{\partial \rho}{\partial x}, \frac{\partial (\rho U)}{\partial x}, \frac{\partial (\rho V)}{\partial x}, \frac{\partial E}{\partial x} \right)^T,
\label{2dsolve_a}
\end{equation}
where
\begin{equation}
\left(\Gamma_{\alpha \beta}\right)=\left(\int g \psi_\alpha \psi_\beta\, \mathrm{d}\Xi/\rho\right)=
 \left(\begin{array}{cccc}
\displaystyle 1 & U &V &\Phi_1\\[3mm]
\displaystyle U & U^2+\frac{1}{2 \lambda} & UV &\Phi_2 \\[3mm]
\displaystyle V & UV & V^2+\frac{1}{2 \lambda}  &\Phi_3 \\[3mm]
\displaystyle \Phi_1 & \Phi_2 & \Phi_3 & \Phi_4
\end{array}\right),
\label{2dGamma}
\end{equation}
with
\begin{equation}
\begin{split}
 & \Phi_1=\frac{1}{2} \left(U^2+V^2+\frac{K+2}{2 \lambda}\right), \quad
 \Phi_2=\frac{1}{2} \left(U^3+UV^2+\frac{(K+4)U}{2 \lambda}\right), \quad \\
 &\Phi_3=\frac{1}{2} \left(U^3+U^2V+\frac{(K+4)V}{2 \lambda}\right), \quad \Phi_4=\frac{1}{4} \left((U^2+V^2)^2+\frac{(K+4)(U^2+V^2)}{\lambda}+\frac{K^2+6K+8}{4 \lambda^2}  \right). \notag
\end{split}
\end{equation}
Thus, $a$ can be obtained by solving the linear system~\eqref{2dsolve_a}. Similar procedures can be used to get $b$.
$A$ is related to the temporal slope {of $g$} with the following form %\QQ{~\cite{xu2001gas}},
\begin{equation}
\frac{\partial_t g}{g} \doteq A=A_1+A_2u+A_3v+A_4 \frac12 (u^2+v^2+|\bm{\xi}|^2), \notag
\end{equation}
where $A_1$, $A_2$, $A_3$, $A_4$ are uniquely determined by the compatibility condition %\QQ{~\cite{xu2001gas}}
\[
\int (au+bv+A)\bm{\psi} g \, \mathrm{d}\Xi=0.
\]
After $a$ and $b$ are determined, $A_1$, $A_2$, $A_3$, $A_4$ can be obtained by solving the following linear system
\begin{equation}
\Gamma (A_1, A_2, A_3, A_4)^T=-\frac{1}{\rho} \int (au+bv) \bm{\psi} g \, \mathrm{d}\Xi,
\label{2dsolve_A}
\end{equation}
with $\Gamma$ specified in equation~\eqref{2dGamma}.

After the distribution function $f$ in equation~\eqref{gasdis2d} is determined, we can get the fluxes ${\mathbf {G}}$ and ${\mathbf {H}}$ in equations~(\ref{semi2d_1}) and (\ref{semi2d_2}) by taking the moments
of $u\bm{\psi}$ and $v \bm{\psi}$, they are
\beq
\begin{split}
&\mathbf{G}=\int u \bm{\psi} f d\Xi=\rho \left[ \langle u\bm{\psi} \rangle-\tau \langle au^2\bm{\psi} \rangle- \tau \langle b u v \bm{\psi} \rangle- \tau \langle Au \bm{\psi} \rangle \right],\\
&\mathbf{H}=\int v \bm{\psi} f d\Xi=\rho \left[ \langle v\bm{\psi} \rangle-\tau \langle a u v \bm{\psi} \rangle- \tau \langle bv^2 \bm{\psi} \rangle- \tau \langle A v \bm{\psi} \rangle \right],
\end{split} \notag
\label{GH2ddeterm}
\eeq
where the evaluation of the moments for the 2D Maxwellian distribution function are organized in Appendix~\ref{APP2}.

\subsection{Numerical boundary conditions}
 \label{BC}

For the numerical tests in this paper, we follow closely the examples in~\cite{xu2001gas, liu2007runge}, with similar boundary conditions, e.g. inflow, outflow and wall boundary conditions.
In the CDG method, two sets of approximate solutions on overlapping cells are updated; thus numerical boundary conditions are needed for both solutions.
The inflow and outflow conditions can be treated, in a similar manner as those in DG, for both solutions in the CDG method.
The more challenging case is the wall boundary conditions.

\begin{figure}[ht]
\begin{center}
\includegraphics[width=3. in]{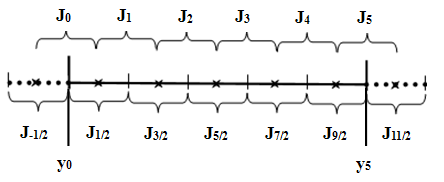}
\end{center} \caption{The 1D overlapping cells for the Couette flow.}
\label{1dcouette}
\end{figure}

In the following, we take the Couette flow in a channel with the bottom wall fixed and the top wall moving (see Section~\ref{coue}) as an example to describe the proposed numerical no-slip boundary condition at both walls. Although the Couette flow is a 2D problem, it can be implemented as a 1D problem since the solutions do not depend on $x$. We assume the overlapping cells in the $y$-direction as plotted in Fig.~\ref{1dcouette} with the walls located at $y_0$ and $y_5$ {with $y_0 = 0$, $y_5 = 5$}. % denote $y=0$ and $y=5$ respectively}.
Cells $J_0=[y_{-\frac{1}{2}}, y_{\frac{1}{2}}]$ and $J_5=[y_{\frac{9}{2}}, y_{\frac{11}{2}}]$ are cut through by the walls.
For the no-slip boundary condition at wall, the physical macroscopic velocities $U, V$ are zero on the wall. However, the numerical ones might not be zero due to numerical errors; such non-zero errors might accumulated during long time evolution and eventually impact the effectiveness of the proposed scheme. We propose to enforce zero velocities at the numerical level.
%\QQ{Suggest to just keep one of examples for $J_0$ or $J_{\frac12}$. There seems no need to include both %examples, since the spirits are similar!}
For example, in cell $J_0$ at the bottom wall, we adopt the following basis functions
\begin{equation}
 \eta_0^0=1, \quad \eta_0^1=\left( \frac{y-y_0}{\Delta y_0 /2} \right), \quad \eta_0^2=\left( \frac{y-y_0}{\Delta y_0 /2} \right)^2, \notag
\end{equation}
with
\begin{equation}
 (\rho U)^Z_h=(\rho U)_0^{Z,0} \eta_0^0+(\rho U)_0^{Z,1} \eta_0^1 +(\rho U)_0^{Z,2} \eta_0^2,\quad (\rho V)^Z_h=(\rho V)_0^{Z,0} \eta_0^0+(\rho V)_0^{Z,1} \eta_0^1 +(\rho V)_0^{Z,2} \eta_0^2, \notag
\end{equation}
where $\Delta y_0=y_{\frac{1}{2}}-y_{-\frac{1}{2}}$, $(\rho U)_0^{Z,l}$ are the coefficients of $(\rho U)^Z_h$
for the basis $\eta_0^l$ with $l=0,1,2$, similarly for $(\rho V)_0^{Z,l}$. Since $\eta_0^1=\eta_0^2=0$ at $y_0$, we only need to enforce $(\rho U)_0^{Z,0}=(\rho V)_0^{Z,0}=0$ to get $(\rho U)_{y_0}=(\rho V)_{y_0}=0$.
Similarly for another set of solutions in cell $J_{\frac{1}{2}}$.
% we adopt the following basis functions
%\begin{equation}
% \zeta_{\frac{1}{2}}^0=1, \quad \zeta_{\frac{1}{2}}^1=\left( \frac{y-y_{\frac{1}{2}}}{\Delta y_{\frac{1}{2}} %/2}+1 \right), \quad \zeta_{\frac{1}{2}}^2=\left( \frac{y-y_{\frac{1}{2}}}{\Delta y_{\frac{1}{2}} /2}+1 %\right)^2, \notag
%\end{equation}
%with
%\begin{equation}
% (\rho U)^W_h=(\rho U)_{\frac{1}{2}}^{W,0} \zeta_{\frac{1}{2}}^0+(\rho U)_{\frac{1}{2}}^{W,1} %\zeta_{\frac{1}{2}}^1 +(\rho U)_{\frac{1}{2}}^{W,2} \zeta_{\frac{1}{2}}^2, \quad
% (\rho V)^W_h=(\rho V)_{\frac{1}{2}}^{W,0} \zeta_{\frac{1}{2}}^0+(\rho V)_{\frac{1}{2}}^{W,1} %\zeta_{\frac{1}{2}}^1 +(\rho V)_{\frac{1}{2}}^{W,2} \zeta_{\frac{1}{2}}^2, \notag
%\end{equation}
%where $\Delta y_{\frac{1}{2}}=y_1-y_0$. $\zeta^{1}_{\f12}=\zeta^{2}_{\f12}=0$ at $y_0$ and we only need to %enforce $(\rho U)_{\frac{1}{2}}^{W,0}=(\rho V)_{\frac{1}{2}}^{W,0}=0$ to get $(\rho U)_{y_0}=(\rho %V)_{y_0}=0$.
The approximate solution in the ghost cell $J_{-\frac{1}{2}}$ is obtained in a mirror-symmetric manner with respect to the solution on cell $J_{\frac{1}{2}}$.
Similar ideas of using a special set of basis to preserve the solution structure in DG methods can also be found in~\cite{cockburn2004locally}.
Boundary conditions on the top wall $y=5$ can be set similarly.

%% file: numerical.tex
\section{Numerical examples}
\label{secnu}
\setcounter{equation}{0}
\setcounter{figure}{0}
\setcounter{table}{0}

In this section, we present simulation results of the proposed CDG-BGK method for several 1D and 2D viscous flow problems.
For comparison, most of the examples are taken from~\cite{xu2001gas, liu2007runge}. The maximum time step
$\Delta \tau^n$ in equations.~(\ref{semi1d_1}), (\ref{semi1d_2}), (\ref{semi2d_1}) and (\ref{semi2d_2}) is chosen based on the CFL condition, while
the time step in evolution is taken to be $\Delta t^n=0.9 \Delta \tau^n$ unless otherwise specified. %{$\Delta t^n=0.5 \Delta \tau^n$ for shock tube problem, $\Delta t^n=0.9 \Delta \tau^n$ for other examples.}
%In order to validate the accuracy and performance of the CDG-BGK scheme, five numerical examples including %1D and 2D problems are tested in this section.
%To show accuracy and performance of the proposed scheme in resolving shock structure and boundary layers,
%Most examples are taken from~\cite{xu2001gas, liu2007runge}.
%The $\Delta \tau^n$ is the maximum time step which is chosen with CFL condition, while $\Delta t^n=0.9 %\Delta \tau^n$ in all examples.

We define, in 1D cases,
\begin{equation}
 \Delta \tau^n=\mbox{min} \left(\frac{\mbox{CFL}_c\, h_x}{\mbox{Em}_x}, \frac{\mbox{CFL}_d\, h_x^2}{\mu} \right), \notag
\end{equation}
and in 2D cases,
\begin{equation}
 \Delta \tau^n=\mbox{min} \left(\mbox{CFL}_c \Big/ \left( \frac{\mbox{Em}_x}{h_x}+ \frac{\mbox{Em}_y}{h_y} \right), \mbox{CFL}_d \Big/ \left( \frac{\mu}{h_x^2}+ \frac{\mu}{h_y^2} \right) \right), \notag
\end{equation}
where
%\QQ{$\mbox{CFL}_c$ and $\mbox{CFL}_d$ are the CFL numbers for convection and diffusion terms respectively,}
 $h_x=\mbox{min}_i \left(\mbox{min}(\Delta x_i, \Delta x_{i+\frac{1}{2}})\right)$ and $h_y=\mbox{min}_j \left(\mbox{min}(\Delta y_j, \Delta y_{j+\frac{1}{2}})\right)$, with $\Delta x_i=x_{i+\frac12}-x_{i-\frac12}$ and $\Delta x_{i+\frac12}=x_{i+1}-x_i$, similarly for $\Delta y_j$ and $\Delta y_{j+\frac12}$. $\mbox{Em}_x, \mbox{Em}_y$ are the maximum eigenvalues in $x$ and $y$ directions respectively.
The eigenvalues are $U-C, U, U, U+C$ in the $x$ direction and $V-C, V, V, V+C$ in the $y$ direction for convection part, where $C=\sqrt{\gamma p/\rho}$ is the speed of sound.
{$\mbox{CFL}_c$ and $\mbox{CFL}_d$ are the CFL numbers for the convection and diffusion parts.
In our numerical examples, for the third-order TVD Runge-Kutta method~\eqref{rk3}, we take $\mbox{CFL}_c=0.58, 0.33, 0.22$ ~\cite{liu20082} and $\mbox{CFL}_d=0.06, 0.01, 0.005$ for $P^1$, $P^2$ and $P^3$ solution spaces respectively. Our $\mbox{CFL}_d$ numbers are larger than those taken in~\cite{liu2011central} for the central local DG method, yet they are working properly for all our numerical examples.}
%\QQ{For the third-order TVD Runge-Kutta method (\ref{rk3}) with $P^1$, $P^2$ and $P^3$ solution spaces, we %take $\mbox{CFL}_c=0.58, 0.33, 0.22$ for $P^1$, $P^2$ and $P^3$ case~\cite{liu20082}, and %$\mbox{CFL}_d=0.06, 0.01, 0.005$ for $P^1$, $P^2$ and $P^3$ case.}
The Prandtl number modification~\eqref{pr_modify} will be used in all numerical examples except the laminar boundary layer case.

%\QQ{For the Couette flow and Laminar boundary layer problem, the solutions are the steady state solution.
%Numerically we need to specify when the solutions have converged to the steady-state. My treatment here is
%a little artificial and problem dependent. Do we need to clearly describe it in our paper? This is for discussion.}
%\QQ{Shall we describe how to judge the convergence of steady state solutions? Xu described steady state %solutions are obtained by a long time calculation in 2001 paper. But how long? I don't know.}

\subsection{Accuracy test}
We first solve the Navier-Stokes equations~\eqref{1d_NS_simplified} with smooth solutions, where the initial
conditions are given by
\begin{equation}
\rho(x,t=0)=1+0.2\sin(\pi x), \quad U(x,t=0)=1,\quad p(x,t=0)=1.
\label{accu_init}
\end{equation}
The computational domain is $[0,2]$ with periodic boundary condition. Two different viscosity coefficients are tested, {$\mu=0.00001$} and $\mu=0.1$, corresponding to a convection-dominated flow and a viscous flow,
respectively. For this example, the Prandtl number is $Pr=2/3$ and the ratio of specific heats is $\gamma=5/3$. We compute the solutions up to time $t=2$. TVB limiter is not used for this
case. Since the exact solution is not available for this problem, the numerical errors and orders of density
$\rho$ are computed by comparing to the reference solution which is obtained by the $P^3$ solution space with 1280 cells. Here we take $\Delta \tau^n=\mbox{min} \left(\frac{\mbox{CFL}_c\, h_x^{\frac{4}{3}}}{\mbox{Em}_x}, \frac{\mbox{CFL}_d\, h_x^2}{\mu} \right)$ for the $P^3$ case so that the temporal error is not dominated. The results are shown in Table~\ref{table1}. $(k+1)^{th}$-order convergent rate can be observed for the proposed CDG-BGK scheme with {$\mu=0.00001$} and $P^k$ solution spaces, while $k^{th}$-order convergent rate for even $k$ and $(k+1)^{th}$-order convergent rate for odd $k$ can be observed for the solution with $\mu=0.1$ and $P^k$ solution spaces.
\begin{table}
\centering
 \caption{Accuracy test, $L^1$ and $L^\infty$ errors and orders for the initial condition~\eqref{accu_init}
 with $P^1$, $P^2$ and $P^3$ solution spaces. }
  \begin{tabular}{|c|c|c|c|c|c|}
    \hline
    &N  &  $L^1$ error & order   & $L^\infty$ error  & order \\\hline
\multirow{5}{*}{$\mu=0.00001, P^1$}
  &10  &0.21E-02  &  --   &0.34E-02  & --  \\  \cline{2-6}
  &20  &0.78E-03  &1.45   &0.12E-02  &1.50 \\ \cline{2-6}
  &40  &0.24E-03  &1.72   &0.36E-03  &1.70 \\ \cline{2-6}
  &80  &0.65E-04  &1.88   &0.10E-03  &1.86 \\ \cline{2-6}
  &160 &0.17E-04  &1.94   &0.26E-04  &1.93  \\  \hline
\multirow{5}{*}{$\mu=0.00001, P^2$}
  &10  &0.28E-03  & --   &0.46E-03  & --  \\  \cline{2-6}
  &20  &0.36E-04  &2.94  &0.58E-04  &2.98 \\ \cline{2-6}
  &40  &0.46E-05  &2.98  &0.73E-05  &3.00 \\ \cline{2-6}
  &80  &0.58E-06  &2.99  &0.92E-06  &3.00 \\ \cline{2-6}
  &160 &0.73E-07  &2.99  &0.12E-06  &2.99  \\  \hline
\multirow{5}{*}{$\mu=0.00001, P^3$}
  &10  &0.27E-04 &  --  &0.38E-04  &  --  \\  \cline{2-6}
  &20  &0.16E-05 &4.04  &0.25E-05  &3.95 \\ \cline{2-6}
  &40  &0.10E-06 &4.02  &0.15E-06  &3.99 \\ \cline{2-6}
  &80  &0.63E-08 &4.01  &0.97E-08  &3.99 \\ \cline{2-6}
  &160 &0.39E-09 &4.00  &0.61E-09  &3.99  \\  \hline
\multirow{5}{*}{$\mu=0.1, P^1$}
  &10  &0.10E-02  &  --   &0.19E-02  & --  \\  \cline{2-6}
  &20  &0.27E-03  &1.93   &0.52E-03  &1.89 \\ \cline{2-6}
  &40  &0.68E-04  &1.98   &0.13E-03  &1.97 \\ \cline{2-6}
  &80  &0.17E-04  &2.00   &0.33E-04  &1.99 \\ \cline{2-6}
  &160 &0.42E-05  &2.02   &0.84E-05  &2.00  \\  \hline
\multirow{5}{*}{$\mu=0.1, P^2$}
  &10  &0.13E-03  & --   &0.22E-03  & --  \\  \cline{2-6}
  &20  &0.35E-04  &1.95  &0.56E-04  &1.94 \\ \cline{2-6}
  &40  &0.87E-05  &1.99  &0.14E-04  &1.99 \\ \cline{2-6}
  &80  &0.22E-05  &2.00  &0.35E-05  &2.00 \\ \cline{2-6}
  &160 &0.55E-06  &2.00  &0.89E-06  &2.00  \\  \hline
\multirow{5}{*}{$\mu=0.1, P^3$}
  &10  &0.93E-05 &  --  &0.17E-04  &  --  \\  \cline{2-6}
  &20  &0.74E-06 &3.65  &0.14E-05  &3.61 \\ \cline{2-6}
  &40  &0.49E-07 &3.92  &0.90E-07  &3.93 \\ \cline{2-6}
  &80  &0.31E-08 &3.99  &0.56E-08  &4.02 \\ \cline{2-6}
  &160 &0.19E-09 &4.03  &0.38E-09  &3.89  \\  \hline
  \end{tabular}
\label{table1}
\end{table}

\subsection{Couette flow}
\label{coue}
In the second example, we consider the couette flow in a channel of height $H$,
with the bottom wall fixed and the top wall moving at a constant speed $U_1$ in
the horizontal direction.
We assume isothermal boundary condition at the bottom and top walls with temperature being $T_0$ and $T_1$ respectively.
%\QQ{How about the following: "We assume isothermal boundary condition at the bottom and top walls with temperature being $T_0$ and $T_1$ respectively." (instead of
% "The temperatures on the bottom and top walls are fixed with values $T_0=1/ \lambda_0$ and $T_1=1/ \lambda_1$ respectively.")}
If the viscosity and heat conduction coefficients
$\mu$ and $\kappa_q$ are constant, an analytical solution for the steady state temperature distribution can be obtained, that is
\begin{equation}
 \frac{T-T_0}{T_1-T_0}=\frac{y}{H}+\frac{PrEc}{2}\frac{y}{H}(1-\frac{y}{H}),
 \label{couette_ana}
\end{equation}
where the Eckert number is $Ec=U_1^2/(C_p(T_1-T_0))$. $C_p$ is the heat capacity at a constant pressure, for a monatomic gas $C_p=\frac52 R$ and for a diatomic gas $C_p=\frac72 R$.

The solution of this problem does not depend on $x$, hence we solve it as a reduced 1D problem from equation~\eqref{2dns} in the $y$ direction, that is,
\begin{equation}
\mathbf{W}_t+\mathbf{H}_y=0, \notag
\end{equation}
i.e. the compressible Navier-Stokes equations~\eqref{2dNavier-Stokes} without the $x$-derivative term.
%\begin{equation}
%\left(\begin{array}{c}
%\displaystyle \rho \\[3mm]
%\displaystyle \rho U\\[3mm]
%\displaystyle \rho V \\[3mm]
%\displaystyle E
%\end{array}\right)_t +
%\left(\begin{array}{c}
%\displaystyle \rho V \\[3mm]
%\displaystyle \rho UV \\[3mm]
%\displaystyle \rho V^2 + p \\[3mm]
%\displaystyle V(E+p) \end{array}\right)_y=
%\left(\begin{array}{c}
%\displaystyle 0 \\[3mm]
%\displaystyle s_{1y} \\[3mm]
%\displaystyle s_{2y} \\[3mm]
%\displaystyle s_{3y} \\[3mm]\end{array}\right)_y, \notag
%\end{equation}
%where $s_{1y}=\mu U_y$, $s_{2y}=\mu \left[ 2 V_y- \frac{2}{K+2} V_y \right]$,
%$s_{3y}=\mu \left[ UU_y+2VV_y- \frac{2}{K+2} VV_y +\frac{K+4}{4} T_y \right]$,
%and it can be written as

We take the computational domain to be $[0,5]$ ($H=5$) and divided by $5$ cells with cell size $\Delta y=1$.
The isothermal no-slip boundary condition with pressure gradient being zero in $y$ direction is used on the bottom and top walls~\cite{xu2001gas}.
Here we consider the temperature $\lambda_0$ and $\lambda_1$ at the boundaries with
different ratios of specific heats $\gamma=5/3,7/5$, different Prandtl numbers $Pr=0.72,1.0$ and different Eckert numbers $Ec=10,50$. For specific settings, see Table~\ref{couette_boundary}.
We take $U_1=0.1$ and $\mu=0.1$. The initial conditions are
\begin{equation}
\rho(y,t=0)=1, \quad U(y,t=0)=0.1, \quad V(y,t=0)=0, \quad M(y,t=0)=0.1, \notag
\end{equation}
where $M=U/C$ is the Mach number.
%From eq.~(\ref{speedsound}), we get $\lambda= \gamma/2$.
%The boundary conditions are in table~\ref{couette_boundary}. Other values can be used in the test.
 \begin{table}
\centering
 \caption{Boundary settings of $\lambda_0$ and $\lambda_1$ for Couette flow.}
  \begin{tabular}{|c|c|c|c|c|}
    \hline
  $Pr$  & $\gamma$ & $Ec$ &$\lambda_0$ &$\lambda_1$    \\ \hline
  0.72,1.0 &5/3 &10  & 1/1.19960  & 1/1.20040   \\  \hline
  0.72,1.0 &5/3 &50  & 1/1.19992  & 1/1.20008  \\  \hline
  0.72,1.0 &7/5 &10  & 1/1.42851  & 1/1.42863    \\ \hline
  0.72,1.0 &7/5 &50  & 1/1.42829  & 1/1.42886    \\  \hline
  \end{tabular}
\label{couette_boundary}
\end{table}

%For a monatomic gas, $\gamma=5/3$, $C_p=5/4$, we set
%\begin{equation}
 %\lambda_0=\begin{cases}
%    1/1.1996,  &\mbox{for} \quad Ec=10,\\
%    1/1.19992, &\mbox{for} \quad Ec=50,
%   \end{cases} \quad
%    \lambda_1=\begin{cases}
%    1/1.2004,  &\mbox{for} \quad Ec=10,\\
%    1/1.20008, &\mbox{for} \quad Ec=50.
%   \end{cases} \notag
%\end{equation}
% For a diatomic gas, $\gamma=7/5$, $C_p=7/4$, we set
%\begin{equation}
 %\lambda_0=\begin{cases}
%    1/(10/7-0.02/70),  &\mbox{for} \quad Ec=10,\\
%    1/(10/7-0.02/350), &\mbox{for} \quad Ec=50,
%   \end{cases} \quad
%    \lambda_1=\begin{cases}
%    1/(10/7+0.02/70),  &\mbox{for} \quad Ec=10,\\
%    1/(10/7+0.02/350), &\mbox{for} \quad Ec=50.
%   \end{cases} \notag
%\end{equation}

\begin{figure}[ht]
\begin{center}
\includegraphics[width=3.in]{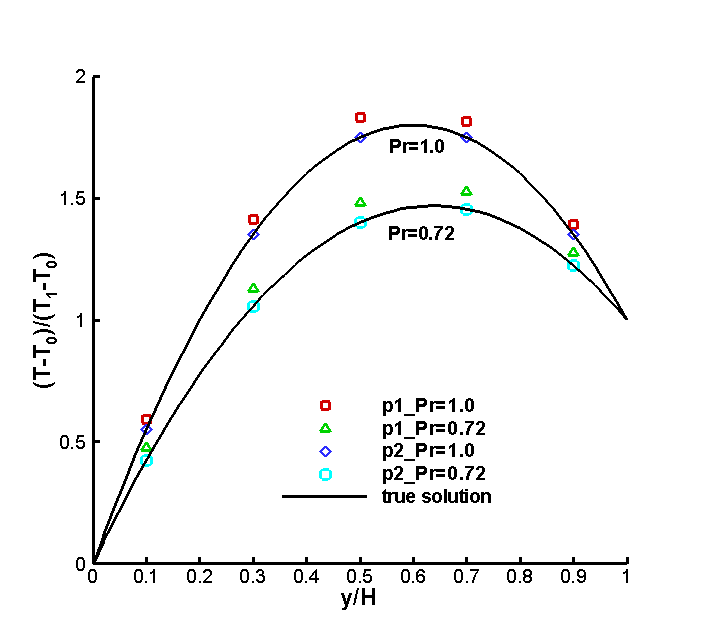}
\includegraphics[width=3.in]{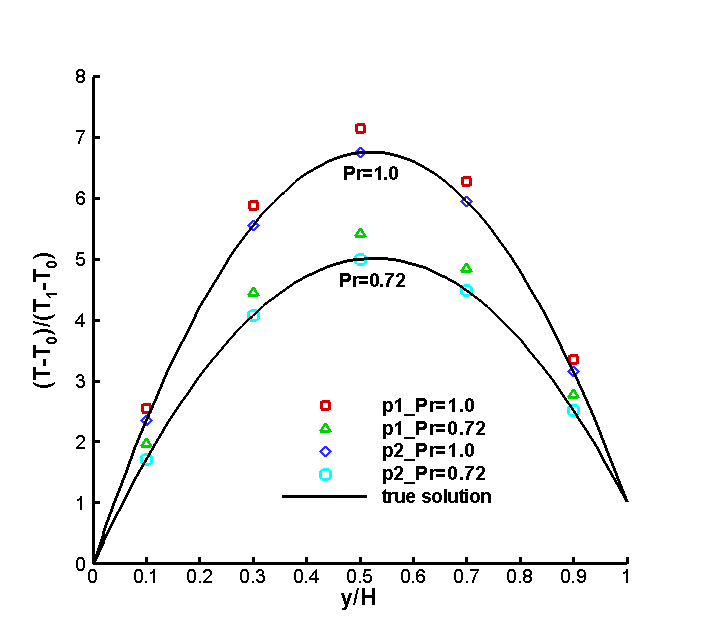}
\end{center} \caption{Temperature ratio $(T-T_0)/(T_1-T_0)$ in the Couette flow with $\gamma=5/3$. Left: $Ec=10$. Right: $Ec=50$.}
\label{figcouette_1}
\end{figure}

\begin{figure}[ht]
\begin{center}
\includegraphics[width=3.in]{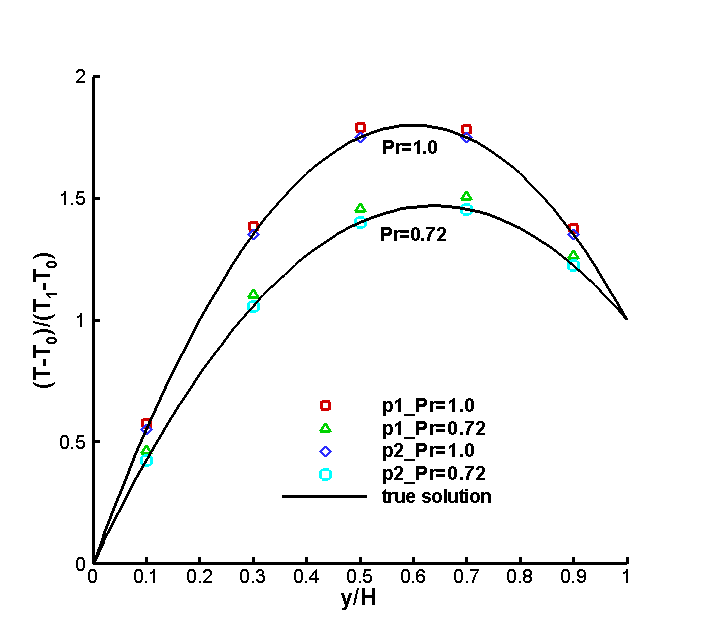}
\includegraphics[width=3.in]{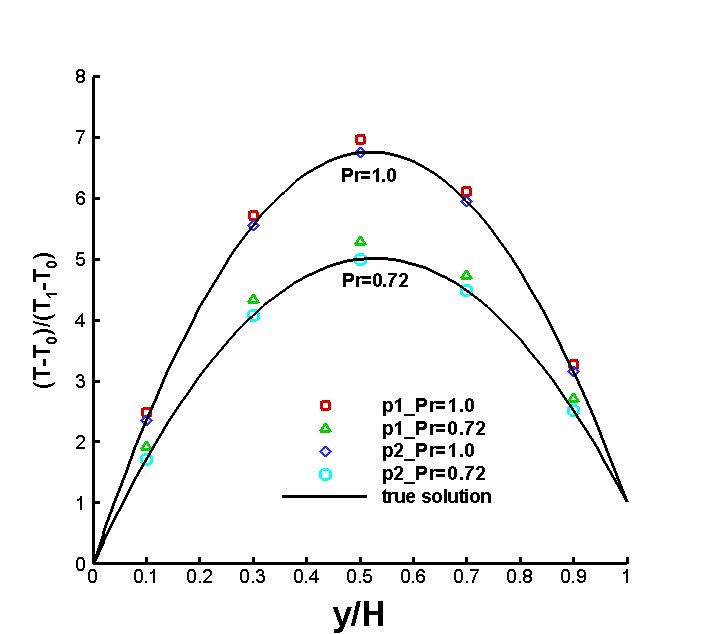}
\end{center} \caption{Temperature ratio $(T-T_0)/(T_1-T_0)$ in the Couette flow with $\gamma=7/5$. Left: $Ec=10$. Right: $Ec=50$.}
\label{figcouette_2}
\end{figure}

The results with different boundary settings are shown in Figs.~\ref{figcouette_1} and \ref{figcouette_2}. The Prandtl number modification~\eqref{pr_modify} is used in $Pr=0.72$ cases.
From the figures, we can see that: (1) numerical results match the analytical solutions very well with different parameters {even on such a coarse numerical cell};
(2) the implementation of the Prandtl number modification is needed compared with the analytical solutions and (3) $P^2$ solution space gives more accurate results than $P^1$ solution space with the same cell size.
Numerical errors and orders of convergence to analytical solutions are summarized in Table~\ref{couette_table}. Roughly $k^{th}$ order convergent rate is observed for the method with $P^k$ polynomial space.
%The order of convergence seems to be sensitive to the boundary treatments.

\begin{table}
\centering
 \caption{Couette flow, $L^1$ and $L^\infty$ errors and orders for $P^1, P^2$ and $P^3$ solution spaces.}
  \begin{tabular}{|c|c|c|c|c|c|}
    \hline
    &N  &  $L^1$ error & order   & $L^\infty$ error  & order \\ \hline
\multirow{4}{*}{$P^1$}
  &5   &0.61E+00  &--    &0.66E+00  &--  \\  \cline{2-6}
  &10  &0.19E+00  &1.66  &0.22E+00  &1.60\\  \cline{2-6}
  &20  &0.61E-01  &1.65  &0.74E-01  &1.55\\  \cline{2-6} \hline
%  &40  &0.22E-01  &1.46  &0.29E-01  &1.36\\  \hline
\multirow{4}{*}{$P^2$}
  &5   &0.13E-03  &--    &0.16E-03  &--  \\  \cline{2-6}
  &10  &0.30E-04  &2.14  &0.41E-04  &1.96\\  \cline{2-6}
  &20  &0.74E-05  &2.04  &0.11E-04  &1.96\\  \cline{2-6} \hline
%  &40  &0.19E-05  &1.99  &0.27E-05  &1.95\\  \hline
\multirow{4}{*}{$P^3$}
  &5   &0.22E-04  &--    &0.28E-04  &--  \\  \cline{2-6}
  &10  &0.35E-05  &2.63  &0.45E-05  &2.62\\  \cline{2-6}
  &20  &0.61E-06  &2.51  &0.82E-06  &2.46\\  \cline{2-6} \hline
%  &40  &0.19E-05  &1.99  &0.27E-05  &1.95\\  \hline
  \end{tabular}
\label{couette_table}
\end{table}

\subsection{Navier-Stokes shock structure}

Now we consider the shock structure problem for a monatomic gas by solving the Navier-Stokes equations~\eqref{1d_NS_simplified}. $\gamma=5/3$ and the dynamical viscosity coefficient is $\mu=\mu_{-\infty} (T/T_{-\infty})^{0.8}$, ${-\infty}$ and $\infty$ denote the
values at the upstream and downsteam respectively. The dynamical viscosity coefficient at the upstream keeps to be a constant $\mu_{-\infty}=0.0005$. The collision time $\tau$ in the BGK model is local via the relationship $\tau = \mu/p$
in each cell. The Mach number $M=1.5$ at the upstream and the Prandtl number $Pr=2/3$. The initial conditions are
\begin{equation}
\left(\begin{array}{c}
\displaystyle \rho \\[3mm]
\displaystyle U \\[3mm]
\displaystyle p
\end{array}\right)_{-\infty} =  \left(\begin{array}{c} \displaystyle 1 \\[3mm]
\displaystyle  1 \\[3mm]
\displaystyle \frac{1}{\gamma M^2}  \end{array}\right), \quad
\left(\begin{array}{c}
\displaystyle \rho \\[3mm]
\displaystyle U \\[3mm]
\displaystyle p
\end{array}\right)_{\infty} =  \left(\begin{array}{c} \displaystyle \frac{(\gamma+1)M^2}{2+(\gamma-1)M^2} \\[3mm]
\displaystyle  \frac{\gamma-1}{\gamma+1}+\frac{2}{(\gamma+1)M^2} \\[3mm]
\displaystyle \left(\frac{2\gamma}{\gamma+1}M^2-\frac{\gamma-1}{\gamma+1}\right) \frac{1}{\gamma M^2}  \end{array}\right). \notag
\end{equation}
%\QQ{Symmetric boundary condition is used on the left and right boundaries. For example, assume %$\mathbf{W}_1, \mathbf{W}_0$ are macroscopic conservative variables on the first cell inside of %computational domain and ghost cell respectively. Then we use
%$\mathbf{W}_0=\mathbf{W}_1$ to implement symmetric boundary condition. We can use it because the shock wave %doesn't arrive the left and right boundaries.  The similar treatment is used to shock tube and laminar %boundary layer cases.}
%old version: %For this problem, at a short time, the shock would have not arrived at the left or right boundary. \QQ{Hence, we can impose the boundary condition as the constant left and right states.
%(suggest; instead of "We take the symmetric boundary conditions on both the left and right boundaries. That is, denoting $\mathbf{W}_1$ to be the macroscopic conservative variables right inside the computational domain adjacent to the boundary and $\mathbf{W}_0$ to be the one on the ghost cell correspondingly, we simply let
%$\mathbf{W}_0=\mathbf{W}_1$. Similar arguments for the following examples with symmetric boundary conditions.")}
For this problem, the shock would not arrive at the left or right boundary.
Hence, we can impose the boundary condition as the constant left and right states.
The reference solution can be obtained by integrating the steady state Navier-Stokes equations, with the corresponding Matlab programs available in Appendix C of~\cite{xu2001gas}.

\begin{figure}[ht]
\begin{center}
\includegraphics[width=3.in]{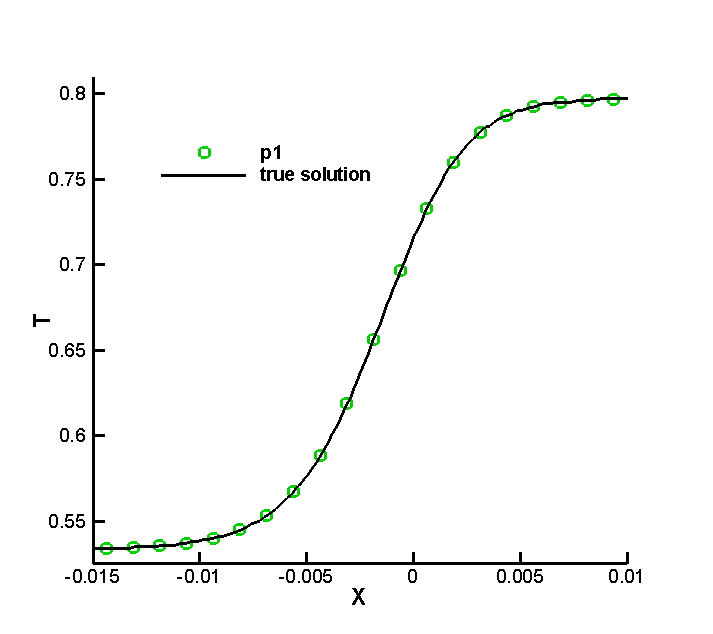}
\includegraphics[width=3.in]{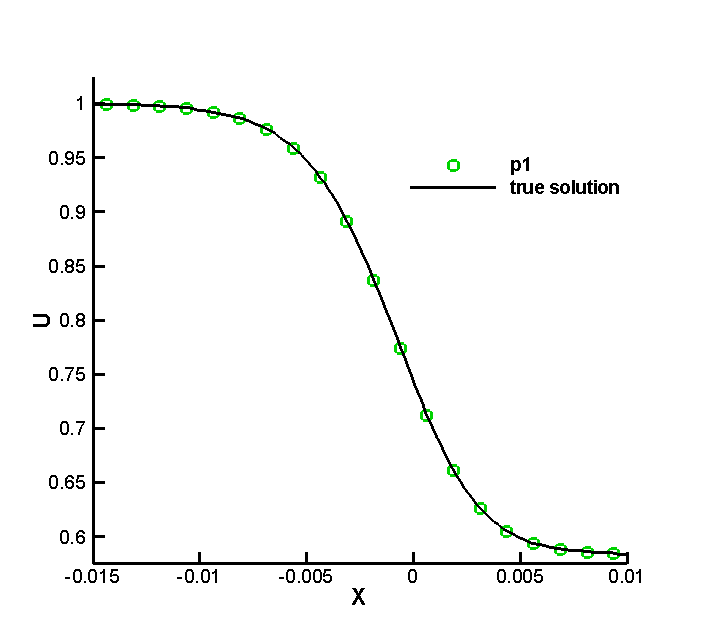}
\includegraphics[width=3.in]{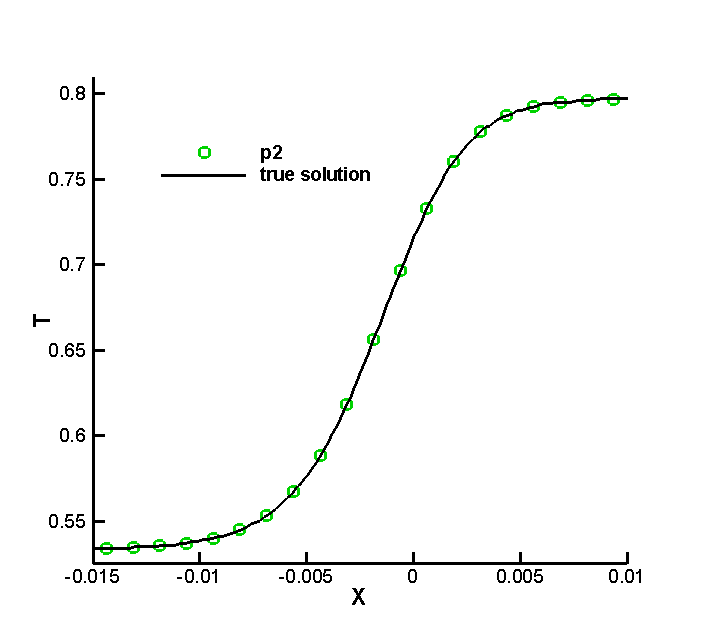}
\includegraphics[width=3.in]{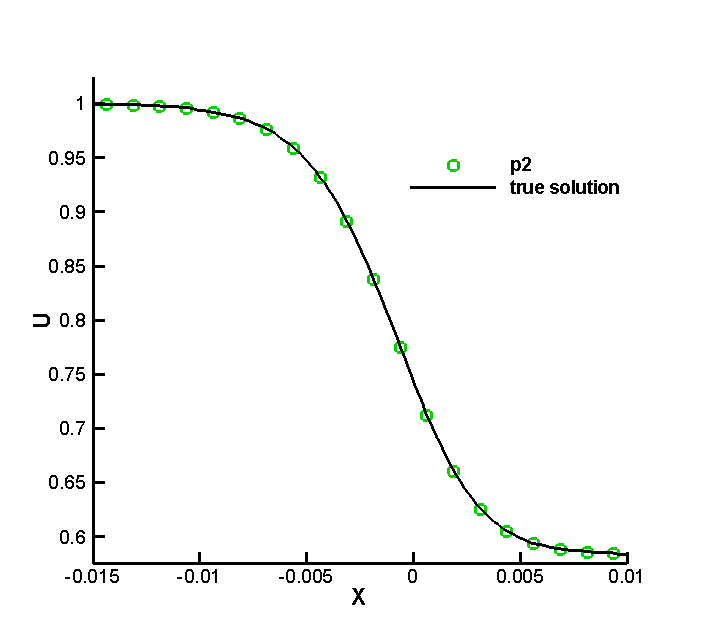}
\includegraphics[width=3.in]{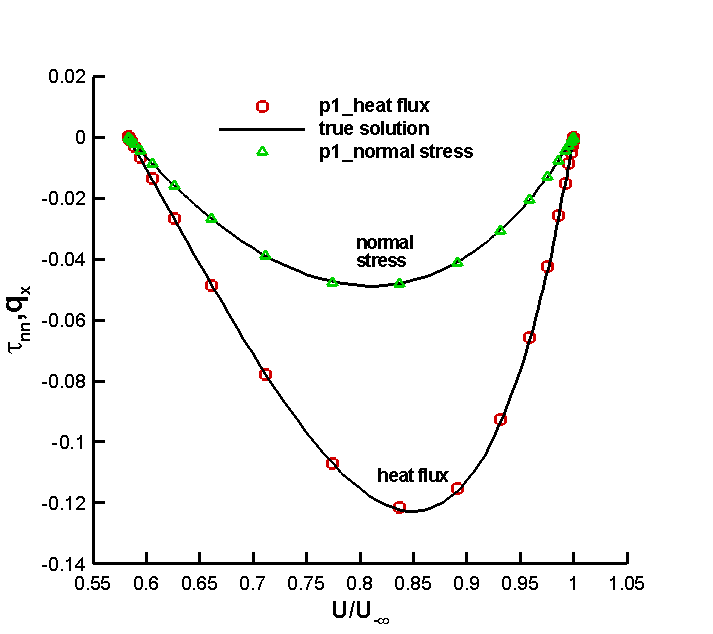}
\includegraphics[width=3.in]{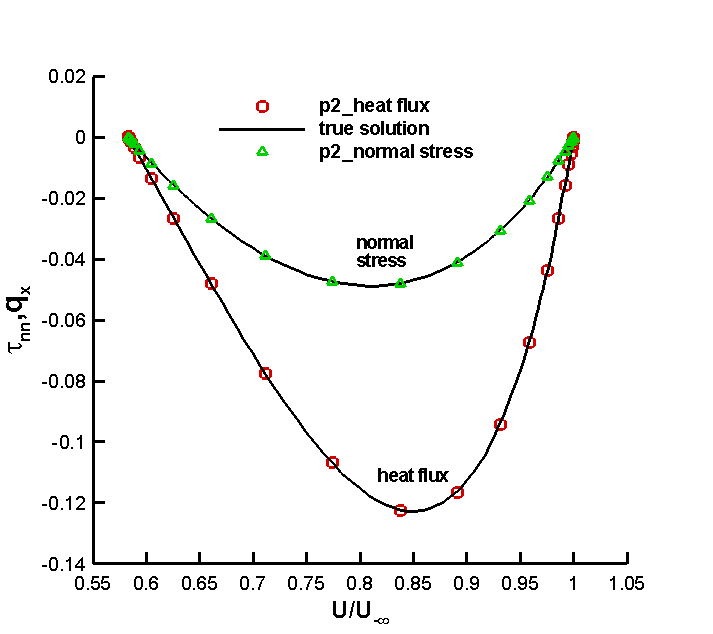}
\end{center} \caption{Navier-Stokes shock structure calculation, $P^1$ and $P^2$ cases.}
\label{shockstructure}
\end{figure}

The computational domain is $[-0.1,0.1]$ and the cell size $\Delta x=1/800$ for both $P^1$ and $P^2$ cases. TVB limiter is used for this example. The results are presented in Fig.~\ref{shockstructure}. In the figures, the normal stress and the heat flux are defined to be
\begin{equation}
\tau_{nn}=\frac{4}{3}\mu\frac{U_x}{2p}, \quad q_x=-\frac{5}{4}\frac{\mu}{Pr}\frac{T_x}{pC}. \notag
\end{equation}
From these results, we can see that the shock structure is captured well with a reasonable number of grid points. The difference between the results from $P^1$ and $P^2$ cases are very small.

\subsection{Shock tube problem}
\begin{figure}[ht]
\begin{center}
\includegraphics[width=3.in]{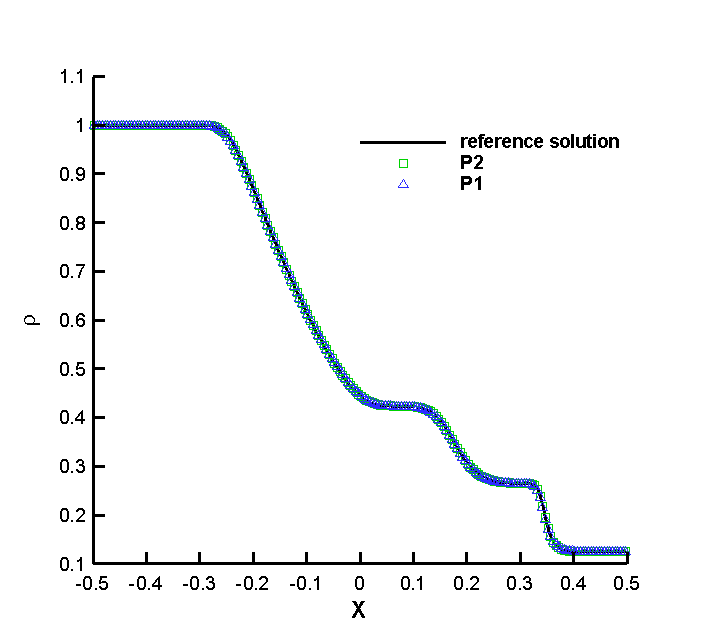}
\includegraphics[width=3.in]{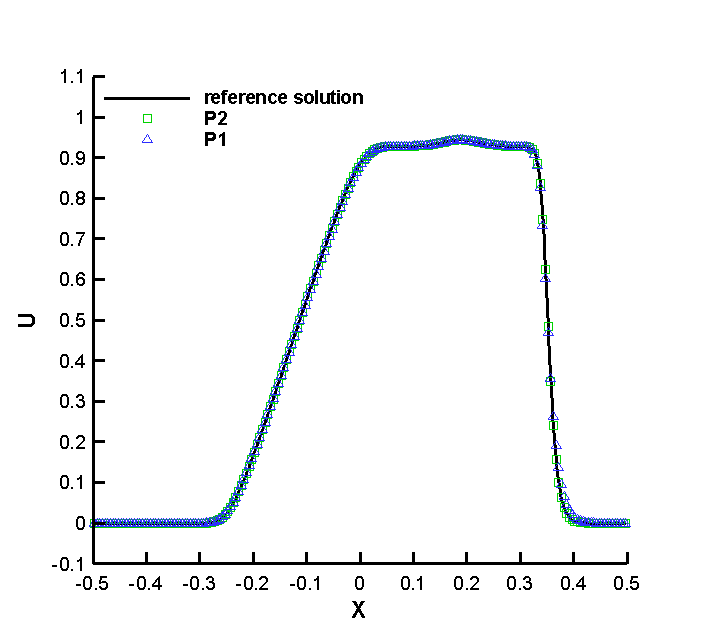} \\
\end{center} \caption{Shock tube problem for the Navier-Stokes equations with kinematic viscosity coefficient $\nu=0.0005/(\rho \sqrt{\lambda})$.}
\label{sod1}
\end{figure}

\begin{figure}[ht]
\begin{center}
\includegraphics[width=3.in]{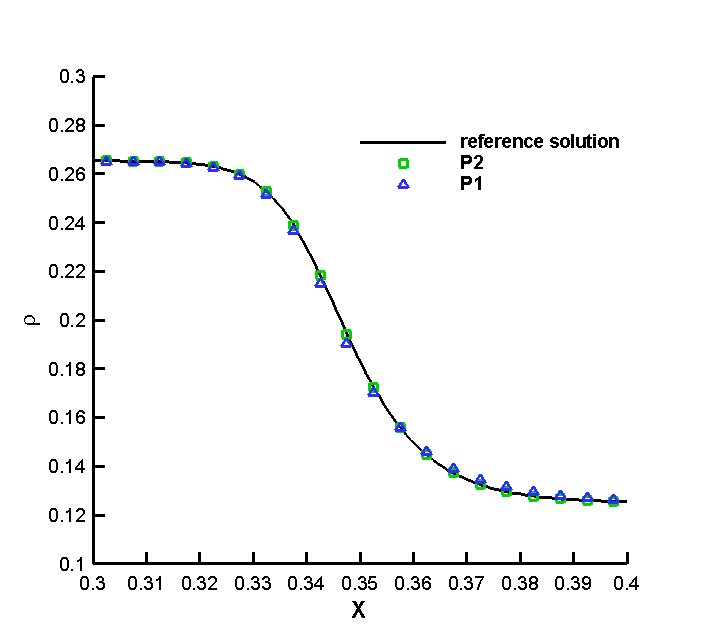}
\end{center} \caption{The zoom-in view of the density distribution around the shock wave in shock tube test with $\nu=0.0005/(\rho \sqrt{\lambda})$.}
\label{sod1refine}
\end{figure}

\begin{figure}[ht]
\begin{center}
\includegraphics[width=3.in]{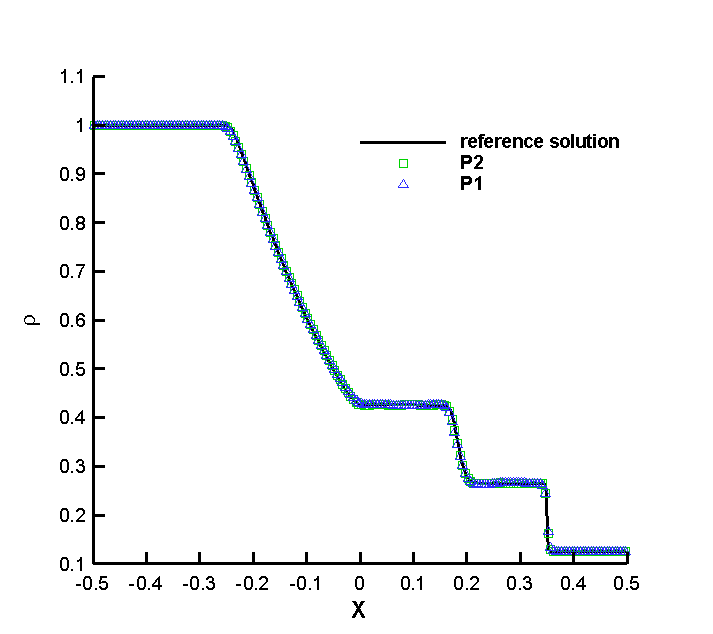}
\includegraphics[width=3.in]{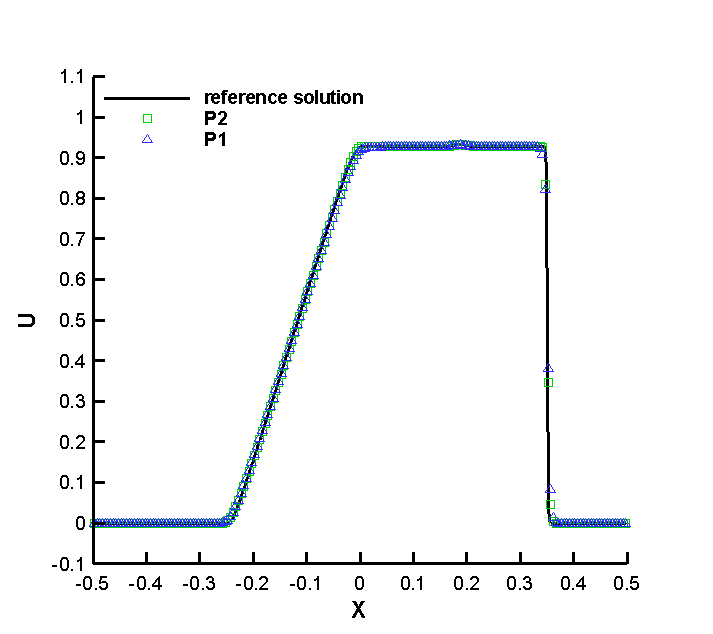} \\
\end{center} \caption{Shock tube problem for the Navier-Stokes equations with kinematic viscosity coefficient $\nu=0.00005/(\rho \sqrt{\lambda})$.}
\label{sod2}
\end{figure}

%In order to further test the CDG-BGK method in capturing the Navier-Stokes solutions in the unsteady case,
In the fourth example, the Sod problem is tested by solving the Navier-Stokes equations~\eqref{1d_NS_simplified} with $\gamma=1.4$ and $Pr=2/3$. The computational domain is $[-0.5,0.5]$ with the
cell size $\Delta x=1/200$, and $\Delta t^n=0.5 \Delta \tau^n$ in this case. TVB limiter is used for this example. The initial conditions are
\begin{equation}
(\rho, U, p)=
\begin{cases}
(1, 0, 1), \quad & x \le 0, \\
(0.125, 0, 0.1), \quad & x > 0.
\end{cases}
\label{sod}
\end{equation}
%Symmetric boundary condition is used on the left and right boundaries.
%\begin{equation}
% \rho=1, \, U=0, \, p=1, \quad \mbox{for} \quad x\le 0; \quad  \rho=0.125, \, U=0, \, p=0.1, \quad %\mbox{for} \quad x\ge 0. \notag
%\end{equation}
We compute the solutions up to time $t=0.2$. 
Similar to shock structure case, we can impose the boundary condition as the constant left and right states.
In Fig.~\ref{sod1}, we show the results with a kinematic
viscosity coefficient $\nu=\mu/ \rho =0.0005/ (\rho \sqrt{\lambda})$.
The solid lines are the reference solutions computed on a much refined cell size $\Delta x=1/1200$ with $P^2$ solution space. Both the shock and the contact discontinuity are captured well. From the zoom-in Fig.~\ref{sod1refine},
we can see the $P^2$ case gives slightly better results than the $P^1$ case.
The results with a smaller viscosity coefficient $\nu =\mu/ \rho =0.00005/ (\rho \sqrt{\lambda})$ are presented in Fig.~\ref{sod2}. %the CDG-BGK scheme becomes a shock capturing scheme (?).
The CDG-BGK method can capture the sharp discontinuity.

\subsection{Laminar boundary layer}
The last example is the 2D laminar boundary layer problem over a flat plate and we solve the 2D Navier-Stokes equations~\eqref{2dns}.
The wall starts from $x=0$ with  a length of $L=3$ at $y=0$. A uniform rectangular mesh with $480 \times 120$ cells is used on the computational domain $[-1,3] \times [0,1]$.
The initial conditions are set to be
\begin{equation}
\rho_{-\infty}=\rho(x,y,t=0)=1, \quad U_{-\infty}=U(x,y,t=0)=1, \quad V(x,y,t=0)=0,\quad M(x,y,t=0)=0.2. \notag
\end{equation}
We take $\gamma=1.4$, $Pr=1$, and $\lambda_{-\infty}=\lambda(x,y,t=0)=\gamma M^2/2$.
The kinematic viscosity coefficient is $\nu=3 \times 10^{-4}$. The Reynolds number based on the upstream flow states and the length $L$ is $Re=\frac{ L U_{-\infty}}{\nu}= 10^4$.
{No limiter is used in this case.} The no-slip adiabatic boundary condition is imposed on the flat plate.
Mirror-symmetric boundary condition is used for the other part of the bottom boundary.
%\QQ{The non-reflecting boundary condition based on the one-dimensional Riemann invariants is used on the %left and top boundaries. For example, on the left boundary, assume ($\rho_1, U_1, p_1, C_1, S_1; \rho_0, %U_0, p_0, C_0, S_0$) are the macroscopic density, velocity in $x$ direction, pressure, speed of sound, %entropy on the first cell inside of computational domain and ghost cell respectively. $C_{-\infty}$ is the %initial speed of sound.
%Based on the Riemann invariants
%\beq
%R_1=U_{-\infty}+\frac{2C_{-\infty}}{\gamma-1},\quad R_2=U_1-\frac{2C_1}{\gamma-1},
%\eeq
%we can determine $U_0=\frac12 (R1+R2)$, $C_0=\frac{\gamma-1}{4} (R1-R2)$. Based on the entropy invariant %condition and the speed of sound relation, we have $S_0=S_1=\frac{p_1}{\rho_1^\gamma}$, %$\rho_0=(\frac{C_0^2}{\gamma \, S_0})^{\frac{1}{\gamma-1}}$, $p_0=\frac{C_0^2 \rho_0}{\gamma}$.
%Finally, we get the macroscopic piecewise constant values ($\rho_0, U_0, p_0$) on ghost cell. Similar %treatment can be used on top boundary.}
On the left and top boundaries, the non-reflective boundary condition is used, which is based on the Riemann invariants.
%\QQ{suggested new version: On the left and top boundaries, the non-reflective boundary condition is used, which is based on the Riemann invariants.}
%\RR{old version: On the left and top boundaries, the reflected wave would not hit these two parts. The non-reflective boundary condition is used, which is based on the Riemann invariants.}
For example, on the left boundary, let $(\rho_1, U_1, p_1)$ be the macroscopic density, velocity and pressure right inside the computation domain adjacent to the left boundary and $(\rho_0, U_0, p_0)$ be the values on the ghost cell correspondingly,
we have
\beq
\rho_0=(\frac{C_0^2}{\gamma \, S_0})^{\frac{1}{\gamma-1}},\quad U_0=\frac12 (R_1+R_2), \quad p_0=\frac{C_0^2 \rho_0}{\gamma}, 
\label{nonreflect}
\eeq
where $R_1=U_{-\infty}+\frac{2C_{-\infty}}{\gamma-1}$ and $R_2=U_1-\frac{2C_1}{\gamma-1}$ are two Riemann invariants, $S_{-\infty}=\frac{p_{-\infty}}{\rho_{-\infty}^\gamma}$ is the entropy from the initial condition, $C_1=\sqrt{\gamma/(2\lambda_1)}$
and $C_{-\infty}=\sqrt{\gamma / (2\lambda_{-\infty})}$ are speeds of sound. $U_0=\frac12 (R_1+R_2), C_0=\frac{\gamma-1}{4} (R_1-R_2), S_0=\frac{p_0}{\rho_0^\gamma}$ in
\eqref{nonreflect} are obtained from the conditions $U_0+\frac{2C_0}{\gamma-1}=R_1$, $U_0-\frac{2C_0}{\gamma-1}=R_2$ and $S_0=S_{-\infty}$.
For details, see~\cite{xu2001dissipative}. Similar treatments on the top boundary.
A first order extrapolation of cell average of the conservative variables is used at the right boundary.

\begin{figure}[ht]
\begin{center}
\includegraphics[width=3.in]{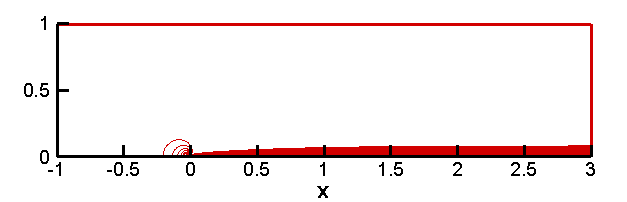}
\includegraphics[width=2.6 in]{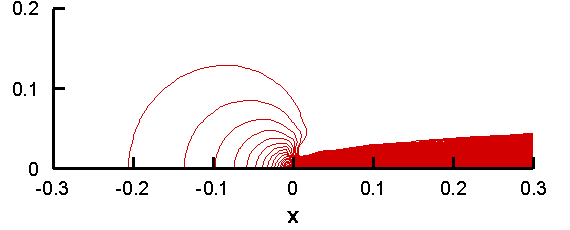}
\end{center} \caption{Laminar boundary layer, contour of velocity obtained by $P^1$ case. The right plot is the zoom-in plot of the interesting region in the left plot.}
\label{laminar_contour}
\end{figure}

\begin{figure}[ht]
\begin{center}
\includegraphics[width=3.in]{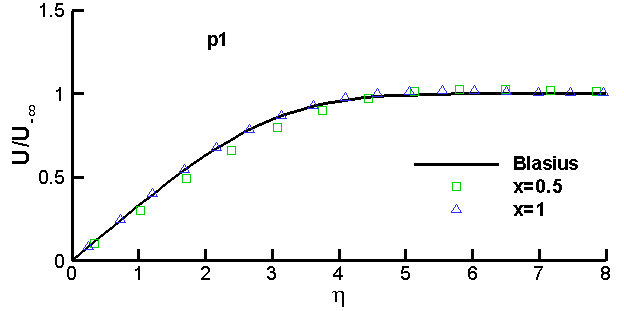}
\includegraphics[width=3.in]{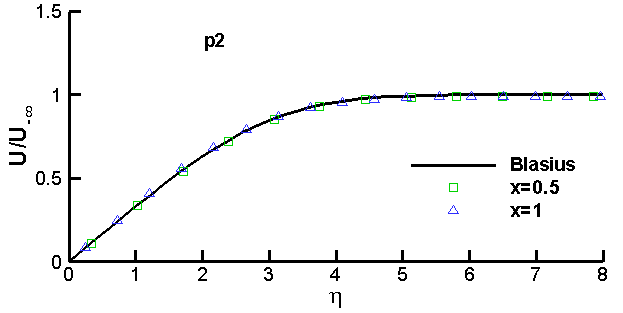}
\end{center} \caption{Laminar boundary layer, $U$ velocity distribution along two vertical lines benchmarked with the Blasius solution. CDG solutions with $P^1$ and $P^2$ solution spaces.}
\label{laminar_nondim}
\end{figure}

In Fig.~\ref{laminar_contour}, we show the $P^1$ solution of the velocity $U$ in the $x$ direction. The non-dimensional velocity $U/U_{-\infty}$ at location $x=0.5$ and $x=1$ for both $P^1$ and $P^2$ cases
are shown in Fig.~\ref{laminar_nondim}, which are compared to the Blasius solution. In the plots, $\eta=y\sqrt{U_{-\infty}/(\nu x)}$, and we can see that the scheme with the $P^2$ solution space performs slightly better than that with  the $P^1$ solution space.

%% file: conclusion.tex
\section{Conclusion}
\label{seccon}
\setcounter{equation}{0}
\setcounter{figure}{0}
\setcounter{table}{0}

In this paper, a novel CDG-BGK method for viscous flow simulations is proposed.
The new scheme inherits several merits from both the CDG framework and the gas-kinetic BGK schemes.
The fluxes in the BGK method is based on the particle transport and collisional mechanism via the gas-kinetic BGK model. Such fluxes take into account of both the convective and viscous terms,
due to the intrinsic connection between the gas-kinetic BGK model and the Navier-Stokes equations. 
The CDG method evolves two pieces of approximate solutions defined on overlapping meshes. The cell interfaces of one computational mesh are inside the staggered mesh, hence the fluxes are in the continuous region of the
staggered solution. For the CDG-BGK method, the distribution function in the interior of elements is continuous and is much easier to evaluate than existing finite volume or DG BGK methods. 
%A third order TVD RK method is applied for high order time discretization. Strategies in handling numerical boundary conditions are proposed and described.
%The costs of obtaining the fluxes is less than will
%comments on CDG: advantages compared with RKDG; and the cost.
Numerical results in 1D and 2D illustrate the accuracy and robustness of the proposed CDG-BGK scheme.

%% file: acknowledge.tex
\bigskip
\noindent
{\bf Acknowledgement.}
Research of the first author is partly supported by the scholarship from China Scholarship Council (CSC) under the Grant CSC No. 201206030032.
The first, third and fourth authors are partially supported by Air Force Office of Scientific Computing YIP grant FA9550-12-0318, NSF grant DMS-1217008 and University of Houston.
The authors would like to thank F.-Y. Li, H.-W. Liu and C.-W. Shu for helpful discussions about boundary conditions.

%% file: appendix.tex
\appendix
\section{Appendix}

\label{appendix}

\renewcommand{\theequation}{A.\arabic{equation}}

%\subsection{Physics notations}
%\label{APP1}
%We describe physical parameters which are used in this paper. For the monatomic gas, the internal degree of freedom $N=0$. For the diatomic gases, $N=2$ accounts for two independent rotational degrees of freedom. Equipartition principle in statistical mechanics
%shows that each degree of freedom shares an equal amount of energy $\frac{1}{2} k_B T$, where $k_B$ is Boltzmann constant. Then the heat capacity of $C_p$ at constant pressure and $C_v$ at constant volume for gases in equilibrium state have the forms
%\begin{equation}
% C_v=\frac{N+3}{2} R; \quad C_p=\frac{(N+3)+2}{2} R, \notag
%\end{equation}
%where $R=k_B/m$ is the gas constant, the 3 accounts for the degree of freedom of molecular motion in $x,y,z$ directions. From the above equations, we can obtain the ratio of specific heats,
%\begin{equation}
% \gamma=\frac{C_p}{C_v}=\frac{(N+3)+2}{N+3}. \notag
%\end{equation}
%For a monatomic gas, $N=0$, $\gamma=5/3$. For a diatomic gas, $N=2$, $\gamma=7/5$.
%The Prandtl number is $Pr=\mu C_p/ \kappa_q=1$ for the BGK model.
%For a monatomic gas,
%the heat conduction coefficient becomes $\kappa_q=5 \mu R/2$.

\subsection{1D and 2D moments}
\label{APP2}
The evaluation of the Maxwellian is given in this section, the details can be found in~\cite{xu2001gas}.
For the 1D flow, the moments of Maxwellian $g$ with respect to $Q$ is introduced as,
\begin{equation}
 \rho \langle Q \rangle=\int Q g du d\bm{\xi}, \notag
\end{equation}
and the general moment formula is
\begin{equation}
 \langle u^n \bm{\xi}^l \rangle=\langle u^n \rangle \langle \bm{\xi}^l \rangle, \notag
\end{equation}
where $n$ is integer, and $l$ is an even integer.
The moments of $\langle \bm{\xi}^l \rangle$ are
\begin{equation}
 \langle \bm{\xi}^0 \rangle=1, \quad
 \langle \bm{\xi}^2 \rangle=\frac{K}{2 \lambda}, \quad
 \langle \bm{\xi}^4 \rangle=\frac{K(K+2)}{4 \lambda^2},
\label{mom_xi}
\end{equation}
and
\begin{equation}
 \langle u^0 \rangle=1, \quad \langle u^1 \rangle=U, \ldots,
 \langle u^{n+2} \rangle=U \langle u^{n+1} \rangle+\frac{n+1}{2 \lambda} \langle u^n \rangle.
 \label{mom_u}
\end{equation}
For the 2D flow,
\begin{equation}
 \rho \langle Q \rangle=\int Q g du dv d\bm{\xi}, \notag
\end{equation}
and the general moment formula is
\begin{equation}
 \langle u^n v^m \bm{\xi}^l \rangle=\langle u^n \rangle \langle v^m \rangle \langle \bm{\xi}^l \rangle, \notag
\end{equation}
where $n,m$ are integers, and $l$ is an even integer. Here the moments $\langle u^n \rangle$ and
$\langle \bm{\xi}^l \rangle$ are the same as the 1D flow. The moments of $v$ are
\begin{equation}
 \langle v^0 \rangle=1, \quad \langle v^1 \rangle=V, \ldots,
 \langle v^{n+2} \rangle=V \langle v^{n+1} \rangle+\frac{n+1}{2 \lambda} \langle v^n \rangle.
 \label{mom_v}
\end{equation}